\documentclass[10pt, a4paper, reqno]{amsart}

\usepackage[english]{babel}
\usepackage{verbatim}
\usepackage[all]{xy}
\usepackage{lunadiagrams}
\usepackage{mathrsfs}
\usepackage{multirow}
\usepackage{hyperref}
\usepackage{amssymb, amsthm, latexsym, mathrsfs, amsfonts, amsxtra, amscd}

\theoremstyle{plain}
\newtheorem*{teo*}{Theorem}
\newtheorem{teo}{Theorem}[section]
\newtheorem{prop}[teo]{Proposition}
\newtheorem{lem}[teo]{Lemma}
\newtheorem{cor}[teo]{Corollary}

\theoremstyle{definition}
\newtheorem{dfn}[teo]{Definition}
\newtheorem{ex}[teo]{Example}

\theoremstyle{remark}
\newtheorem{oss}[teo]{Remark}

\newcommand{\calC}{\mathcal C}
\newcommand{\calL}{\mathcal L}
\newcommand{\calO}{\mathcal O}
\newcommand{\calV}{\mathcal V}
\newcommand{\calX}{\mathcal X} 

\newcommand{\mbA}{\mathbf A}
\newcommand{\mbN}{\mathbf N}
\newcommand{\mbP}{\mathbf P}
\newcommand{\mbQ}{\mathbf Q}
\newcommand{\mbZ}{\mathbf Z}

\newcommand{\scrS}{\mathscr S}

\newcommand{\gra}{\alpha} 
\newcommand{\grb}{\beta}
\newcommand{\grg}{\gamma}
\newcommand{\grd}{\delta}
\newcommand{\grl}{\lambda}     
\newcommand{\grs}{\sigma}
\newcommand{\gro}{\omega}      

\newcommand{\grG}{\Gamma}
\newcommand{\grD}{\Delta}
\newcommand{\grL}{\Lambda} 
\newcommand{\grO}{\Omega}
\newcommand{\grS}{\Sigma}

\newcommand{\sthat}       {\, : \,}
\newcommand{\mand}     {\text{ and }}
\newcommand{\incluso}    {\hookrightarrow}
\newcommand{\lra}        {\longrightarrow}
\newcommand{\vuoto}      {\varnothing}
\renewcommand{\geq}      {\geqslant}
\renewcommand{\leq}      {\leqslant}
\newcommand{\senza}      {\smallsetminus}
\newcommand{\ristretto}  {\bigr|}
\newcommand{\ol}         {\overline}
\newcommand{\wt}         {\widetilde}
\newcommand{\DGH}        {\grD(G/H)}
\newcommand{\Xlt}        {\widetilde{X}_\grd}
\newcommand{\Xl}         {X_\grd}

\newcommand{\BrI}        {{\sf B}^\mathrm{I}_r}
\newcommand{\BrII}       {{\sf B}^\mathrm{I\!I}_r}
\newcommand{\BdueI}      {{\sf B}^\mathrm{I}_2}

\newcommand{\GI}         {{\sf G}^\mathrm{I}_2}
\newcommand{\GII}        {{\sf G}^\mathrm{I\!I}_2}

\newcommand{\quot}[2]{\ensuremath{\raisebox{.40ex}{\ensuremath{#1}} \! \big / \! \raisebox{-.60ex}{\ensuremath{#2}}}}
\newcommand{\Quot}[2]{\ensuremath{\raisebox{.75ex}{\ensuremath{#1}} \! \Big / \! \raisebox{-.75ex}{\ensuremath{#2}}}}

\DeclareMathOperator{\card}{card}
\DeclareMathOperator{\Aut}{Aut}
\DeclareMathOperator{\End}{End}
\DeclareMathOperator{\Hom}{Hom}
\DeclareMathOperator{\Pic}{Pic}
\DeclareMathOperator{\Supp}{Supp}
\DeclareMathOperator{\Stab}{Stab}

\begin{document}

\title[Spherical orbit closures in simple projective spaces]
{Spherical orbit closures in\\simple projective spaces and\\their normalizations}

\author{Jacopo Gandini}

\date{\today}

\curraddr{\textsc{Dipartimento di Matematica ``Guido Castelnuovo''\\
                    ``Sapienza'' Universit\`a di Roma\\
                    Piazzale Aldo Moro 5\\
                    00185 Roma, Italy}}

\email{gandini@mat.uniroma1.it}

\begin{abstract}
Let $G$ be a simply connected semisimple algebraic group over an algebraically closed field $k$ of characteristic 0 and let $V$ be a rational simple $G$-module. If $G/H \subset \mathbf{P}(V)$ is a spherical orbit and if $X = \overline{G/H}$ is its closure, then we describe the orbits of $X$ and those of its normalization $\widetilde{X}$. If moreover the wonderful completion of $G/H$ is strict, then we give necessary and sufficient combinatorial conditions so that the normalization morphism $\widetilde{X} \to X$ is a homeomorphism. Such conditions are trivially fulfilled if $G$ is simply laced or if $H$ is a symmetric subgroup.
\end{abstract}

\maketitle

\section*{Introduction.}

Let $G$ be a simply connected semisimple algebraic group over an algebraically closed field $k$ of characteristic 0; all $G$-modules considered in the following will be supposed to be rational. An algebraic $G$-variety is said to be \textit{spherical} if it is normal and if it contains an open $B$-orbit, where $B\subset G$ is a Borel subgroup; a subgroup $H\subset G$ is said to be \textit{spherical} if the homogeneous space $G/H$ is so: any spherical variety can thus be regarded as an open embedding of a spherical homogeneous space, namely its open $G$-orbit. Important classes of spherical varieties are that of toric varieties and that of symmetric varieties: toric varieties are those spherical varieties whose open orbit is an algebraic torus; symmetric varieties are those spherical varieties whose generic stabilizer $H$ is such that $G^\grs \subset H \subset N_G(G^\grs)$, where $\grs : G \to G$ is an algebraic involution and where $G^\grs$ is the set of its fixed points. Other important classes of spherical varieties are that of flag varieties and the more general one of wonderful varieties: a \textit{wonderful variety} (of \textit{rank} $r$) is a smooth projective $G$-variety having an open $G$-orbit which satisfies following properties:
\begin{itemize}
\item[-] the complement of the open $G$-orbit is the union of $r$ smooth prime divisors having a non-empty transversal intersection;
\item[-] any orbit closure equals the intersection of the prime divisors containing it.
\end{itemize}
A spherical subgroup $H$ is said to be \textit{wonderful} if $G/H$ possesses a wonderful completion (which is unique, if it exists). By \cite{CP} every self-normalizing symmetric subgroup is wonderful; more generally every self-normalizing spherical subgroup is wonderful by \cite{Kn2}.

However many natural examples of embeddings of a spherical homogeneous space do not need to be normal. For instance, consider a simple $G$-module $V$ (in which case we will call $\mbP(V)$ a \textit{simple projective space}) possessing a line $[v]$ fixed by a spherical subgroup. Then consider the orbit $G[v] \subset \mbP(V)$, which is spherical, and take its closure $X = \ol{G[v]} \subset \mbP(V)$, which generally is not normal; denote $\wt{X}$ its normalization. The aim of this work is the study of the orbits of compactifications which arise in such a way, and as well the study of the orbits of their normalizations.

In \cite{BrL} it has been proved that any spherical subgroup which occurs as the stabilizer of a point in a simple projective space is wonderful. If $M$ is the wonderful completion of $G[v]$, then the morphism $G[v] \to X$ extends to $M$ and thus we get a morphism $M\to \wt{X} \to X$: examining such morphism we get a description of the set of orbits of $X$ and of $\wt{X}$. Moreover this leads to a combinatorial criterion to establish whether or not two orbits in $M$ map onto the same orbit in $X$, which in particular implies that different orbits in $X$ are never $G$-equivariantly isomorphic. 

Our main theorem (Theorem \ref{teorema caso stretto}) is a combinatorial criterion for $\wt{X} \to X$ to be bijective; this is done under the assumption that $M$ is \textit{strict}, i.e. that all isotropy groups of $M$ are self-normalizing: strict wonderful varieties, introduced in \cite{Pe}, are those wondeful varieties which can be embedded in a simple projective space; they form an important class of wonderful varieties which generalize the symmetric ones of \cite{CP}. The condition of bijectivity involves the double links of the Dynkin diagram of $G$ and it is trivially fulfilled whenever $G$ is simply laced or $M$ is symmetric; it is easily read off by the \textit{spherical diagram} of $M$, which is a useful tool to represent a wonderful variety starting from the Dynkin diagram of $G$. Main examples of strict wonderful varieties where bijectivity fails arise from the context of \textit{wonderful model varieties} introduced in \cite{L2}; the general strict case is substantially deduced from the model case.

A \textit{model space} for a connected (possibly non-simply connected) semisimple algebraic group $G'$ is a quasi-affine homogeneous space whose coordinate ring contains every simple $G'$-module with multiplicity one; model spaces were classified in \cite{L2}, where it is introduced the \textit{wonderful model variety} $M^{\mathrm{mod}}_{G'}$, whose orbits naturally parametrize up to isomorphism the model spaces for $G'$: every orbit of $M^{\mathrm{mod}}_{G'}$ is of the shape $G'/N_{G'}(H)$, where $G'/H$ is a model space, and conversely this correspondence gives a bijection up to isomorphism.

In order to illustrate the above mentioned criterion of bijectivity in the case of a wonderful model variety, let's set up some further notation.
If $\grl$ is a dominant weight (w.r.t. a fixed maximal torus $T\subset G$ and a fixed Borel subgroup $B \supset T$), define the $\textit{support}$ of $\grl$ as the set $$\Supp(\grl) = \{ \gra \in S \sthat \langle \gra^\vee, \grl \rangle \neq 0\},$$ where $S$ is the set of simple roots w.r.t. $T\subset B$.
If $G_i\subset G$ is a simple factor of type $\sf{B}$ or $\sf{C}$, number the associated subset of simple roots $S_i = \{\alpha^i_1, \ldots \alpha^i_{r(i)}\}$ starting
from the extreme of the Dynkin diagram of $G_i$ which contains the double link; define moreover $S_i^{\mathrm{even}}, S_i^{\mathrm{odd}} \subset S_i$ as the
subsets whose element index is respectively even and odd.
If they are defined, set
\[
            e_i(\lambda) = \min\{k \leqslant r(i) \, : \, \alpha^i_k \in \mathrm{Supp}(\lambda) \cap S_i^{\mathrm{even}}\}
\]
\[
            o_i(\lambda) = \min\{k \leqslant r(i) \, : \, \alpha^i_k \in \mathrm{Supp}(\lambda) \cap S_i^{\mathrm{odd}} \}
\]
or set $e_i(\lambda) = +\infty$ (resp. $o_i(\lambda) = +\infty$) otherwise.
Finally, if $G_i$ is of type $\sf{F}_4$, number the simple roots in $S_i = \{\alpha^i_1, \alpha^i_2, \alpha^i_3, \alpha^i_4\}$ starting from the extreme
of the Dynkin diagram which contains a long root.

\begin{teo*}[see Thm. \ref{teorema caso stretto}]
Suppose that $G[v]\subset \mbP(V)$ is the open orbit of a wonderful model variety $M^{\mathrm{mod}}_{G'}$, where $G'$ is isogenous with $G$; denote $\grl$ the highest weight of $V$ and set $X = \ol{G[v]}$. Then the normalization $\wt{X}\to X$ is bijective if and only if the following conditions are fulfilled, for every connected component $S_i \subset S$:
\begin{itemize}
    \item[i)] If $S_i$ is of type $\sf{B}$, then either $\alpha^i_1 \in \mathrm{Supp}(\lambda)$ or
                $\mathrm{Supp}(\lambda) \cap S^{\mathrm{even}}_i = \varnothing$;
    \item[ii)] If $S_i$ is of type $\sf{C}$, then $o_i(\lambda)  \geqslant e_i(\lambda) - 1$;
    \item[iii)] If $S_i$ is of type $\sf{F}_4$ and $\alpha^i_2 \in \mathrm{Supp}(\lambda)$, then $\alpha^i_3 \in \mathrm{Supp}(\lambda)$ as well.
\end{itemize}
\end{teo*}

When the generic stabilizer $H$ is a self-normalizing symmetric subgroup, compactifications in simple projective spaces were studied in \cite{M}. Under this assumption, an explicit description of the orbits of $X$ was given and it was proved that these orbits are equal to those of the normalization of $X$. Thus our results generalize those contained in \cite{M}.

In the case of a compactification of the adjoint group $G_{\mathrm{ad}}$ (regarded as a $G\times G$-symmetric variety) obtained as the closure of the orbit of the line generated by the identity in a projective space $\mbP(\End(V))$ (where $V$ is a simple $G$-module), a complete classification of the normality and of the smoothness has been given in \cite{BGMR}.

The paper is organized as follows. In section 1, we set notations and preliminaries; in section 2 we give some general results about spherical orbit closures in projective spaces; in section 3 we recall some results from \cite{BrL} about stabilizers of points in simple projective spaces. In section 4, we describe the orbits of the compactifications $X$ and $\wt{X}$; in section 5, we prove the criterion of bijectivity of the normalization map in the strict case; in section 6, we briefly consider the non-strict case giving some sufficient conditions of bijectivity and non-bijectivity of the normalization map.\\

\textit{Aknowledgements.} I want to thank A. Maffei, who proposed me the problem, for all his precious help, and P. Bravi for many useful discussions on the subject. As well, I want to thank the referees for their careful reading and useful comments.

Spherical diagrams have been made with the package \textit{lunadiagrams}, made by P. Bravi and available at http://www.mat.uniroma1.it/$\sim$bravi/lunadiagrams.

\section{Preliminaries.}

Fix a Borel subgroup $B\subset G$ and a maximal torus $T\subset B$; denote $\Phi$ the corresponding root system and $S\subset \Phi$ the corresponding set of simple roots. If $H \subset G$ is any subgroup, denote $\calX(H)$ its character group; if $V$ is a $G$-module, denote $V^{(H)}$ the set of $H$-eigenvectors of $V$ and, if $\chi \in \calX(H)$, denote $V^{(H)}_\chi$ the subset of $V^{(H)}$ where $H$ acts by $\chi$. If $\grl\in \calX(B)$ is a dominant weight, we will denote $V_\grl$ the simple $G$-module with highest weight $\grl$. If $\grL$ is a lattice (i.e. a finitely generated free $\mbZ$-module), then $\grL^\vee=\Hom_\mbZ(\grL,\mbZ)$  denotes the dual lattice and $\grL_\mbQ = \grL \otimes \mbQ$ denotes the rational vector space generated by $\grL$. If $\calC$ is a cone contained in some vector space $V$, then $\calC^\vee$ denotes the dual cone in the dual vector space $V^*$.

Let $X$ be a spherical $G$-variety and fix a base point $x_0\in X$ in such a way that $Bx_0 \subset X$ is an open subset; denote $H = \Stab(x_0)$. Let's introduce some data associated to $X$:
\begin{itemize}

    \item[(1)] $\grL_{X} = \left\{ \textrm{$B$-weights of rational $B$-eigenfunctions in $k(X)$} \right\} \simeq k(X)^{(B)}/k^*$.

    \item[(2)] $\grD(X) = \left\{ \textrm{$B$-stable prime divisors in $X$ which are not $G$-stable} \right\}$, its elements are called the \textit{colors} of $X$. If $Y\subset X$ is an orbit, then $\grD_Y(X)$ denotes the set of colors which contain $Y$.

    \item[(3)] If $\nu : k(X)^* \to \mbQ$ is a rational discrete valuation of $k(X)$, then $\nu$ defines an element $\rho_X(\nu) \in (\grL_{X}^\vee)_\mbQ$ by $$\langle \rho_X(\nu), \chi \rangle = \nu(f_\chi),$$ where $f_\chi \in k(X)^{(B)}$ is any $B$-semiinvariant function of weight $\chi$: since $X$ possesses an open $B$-orbit, such definition does not depend on the function, but only on the weight. If $D\in \grD(X)$, by abuse of notation we will denote $\rho_X(D) = \rho_X(\nu_D)$ the image of the respective valuation $\nu_D$.
    
    A rational discrete valuation $\nu$ is said $G$-invariant if $\nu(g f) = \nu(f)$, for any $f\in k(X)$ and for any $g\in G$; set $\calV_{X}$ the set of $G$-invariant rational valuations of $k(X)$. The map $\rho_X : \calV_X \to (\grL_{X}^\vee)_\mbQ$ identifies $\calV_X$ with a convex cone which generates $(\grL_{X}^\vee)_\mbQ$ as a vector space \cite[Cor. 4.1]{BP}; together with such embedding, $\calV_{X}$ is called the $G$-\textit{invariant valuation cone of} $X$.

\end{itemize}

Both $\grL_X$ and $\grD(X)$, as well as the map $\rho_X$, depend only on the open orbit $G/H \subset X$ and they are the main objects of the \textit{Luna-Vust Theory} (see \cite{Kn}), which classifies normal equivariant embeddings of a given spherical homogeneous space.
A spherical variety is said to be \textit{simple} if it possesses only one closed orbit; it is said to be \textit{toroidal} if no color contains a closed orbit. If a spherical homogeneous space $G/H$ possesses a complete, simple and toroidal embedding, then this is uniquely determined and it is called the \textit{canonical embedding} of $G/H$; we will denote it $M(G/H)$ and it dominates any simple complete embedding of $G/H$. In general, a canonical embedding of $G/H$ exists if and only if the index of $H$ in its normalizer is finite, in which case $H$ is said to be \textit{sober}.

If $X$ is a simple spherical variety with complete closed orbit $Y$, then the Picard group $\Pic(X)$ is freely generated by the classes $[D]$, with $D \in \grD(X)\senza \grD_Y(X)$; moreover, a divisor is generated by  global sections (resp. ample) if and only if it is equivalent to a linear combination of such colors with non-negative (resp. positive) coefficients \cite[Prop. 2.6 and Thm. 2.6]{B1}.

Wonderful varieties are always spherical (see \cite{L0}) and a spherical variety is wonderful if and only if it is complete, toroidal, simple  and smooth. A spherical subgroup which appears as the generic stabilizer of a wonderful variety is said \textit{wonderful}.

If $H$ is a spherical subgroup, then the normalizer $N_G(H)$ acts on the right on $G/H$ by $n\cdot gH = gn^{-1}H$. Consider the induced action of $N_G(H)$ on $\grD(G/H)$: the kernel of such action is called the \textit{spherical closure} of $H$; if $H$ coincides with its spherical closure, then it is called \textit{spherically closed}. Spherically closed subgroups are always wonderful \cite[Cor. 7.6]{Kn2}; a wonderful variety is said to be \textit{spherically closed} if its generic stabilizer is so.

Suppose now $M$ is a wonderful variety with open $B$-orbit $Bx_0$ and set $H = \Stab(x_0)$; suppose moreover that the center of $G$ acts trivially on $M$. Denote $\grD$ the set of colors and $Y\subset M$ the closed orbit; denote $z\in Y^{B^-}$ the $B^-$-fixed point (where $B^-$ denotes the opposite Borel subgroup of $B$ with respect to $T$). Since $G$ is semisimple and simply connected, $\Pic(Y)$ is identified with a sublattice of $\calX(B)$, while $\Pic(G/H)$ is identified with $\calX(H)$: if $\calL \in \Pic(Y)$, then $\calL$ will be identified with the character of $B^-$ acting on the fiber of $\calL$ over $z$, while if $\calL \in \Pic(G/H)$, then $\calL$ will be identified with the character of $H$ acting on the fiber over $eH$.

Let's introduce some more data attached to a wonderful variety $M$, together with some results which can be found with more details and references in \cite{L1} and in \cite{BrL}.
\begin{itemize}

    \item[(4)] $\grS = \left\{ \textrm{$T$-weights of the $T$-module } T_z M/T_z Y \right\}$; its elements are called the \textit{spherical roots} of $M$. Spherical roots form a basis of the lattice $\grL_{G/H}$; they also coincide with the minimal set of generators of the free semigroup $\grL_{G/H}\cap -\calV^\vee_{G/H}$. The cardinality of $\grS$ coincides with the \textit{rank} of $M$, i.e. with the number of $G$-stable prime divisors of $M$, which are naturally in correspondence with spherical roots.  If $\grs\in \grS$, denote $M^\grs$ the corresponding $G$-stable prime divisor: it is a wonderful $G$-subvariety whose set of spherical roots is $\grS \senza \{\grs\}$. Spherical roots are either positive roots or a sum of two positive roots; the \textit{support} of a spherical root $\grs$ is the set of simple roots where $\grs$ is supported.

    \item[(5)]  The \textit{Cartan pairing} of $M$ is the natural pairing $c : \grD \times \grS \to \mbZ$ between colors and spherical roots defined by the equality $[M^\grs] = \sum_{D \in \grD} c(D,\grs) [D]$ in the Picard group $\Pic(M) = \mbZ\grD$. If $D \in \grD$ and $\grs \in \grS$, then $c(D,\grs) = \langle \rho_{G/H}(D), \grs \rangle$.

    \item[(6)] $\grD(\gra)$ = $\left\{ D \in \grD \sthat P_\gra D \neq D \right\}$ is the set of colors \textit{moved by} $\gra$, where $\gra \in S$ and where $P_\gra\supset B$ is the minimal parabolic subgroup associated to $\gra$. For every $\gra\in S$ it holds $0 \leq \card \grD(\gra) \leq 2$.

    \item[(7)] $S^p = \{\gra \in S \sthat \grD(\gra) = \vuoto\}$; it coincides with the set of simple roots associated to the stabilizer of the open $B$-orbit, which is a parabolic subgroup. As well, $S^p$ coincides with the set of simple roots associated with the stabilizer of the $B^-$-fixed point $z$ in the closed orbit $Y$.

    \item[(8)] $S^a = \left\{ \gra \in S \sthat \card \grD(\gra) = 2 \right\} = S \cap \grS$ is the set of simple spherical roots; correspondingly, set $\mbA = \bigcup_{S^a} \grD(\gra)$ the \textit{set of colors of type a}. If $\gra \in S^a$, set $\grD(\gra) = \left\{ D_\gra^+ , D_\gra ^- \right\}$; then for every spherical root $\grs$ it holds
        \[
            c(D_\gra^+, \grs ) + c(D_\gra^-, \grs ) = \langle \gra^\vee , \grs \rangle.
        \]

    \item[(9)] $S^{2a} = \left\{ \gra \in S \sthat 2\gra \in \grS \right\}$; correspondingly, set $\grD^{2a} = \bigcup_{S^{2a}} \grD(\gra)$ the \textit{set of colors of type 2a}. If $\gra \in S^{2a}$, set $\grD(\gra) = \left\{ D_\gra \right\}$; then for every spherical root $\grs$ it holds
        \[
            c(D_\gra, \grs ) = \langle \gra^\vee , \grs \rangle / 2.
        \]

    \item[(10)] $S^b = S \senza (S^p \cup S^a \cup S^{2a})$; correspondingly, set $\grD^b = \bigcup_{S^b} \grD(\gra)$ the \textit{set of colors of type b}. If $\gra \in S^b$, set $\grD(\gra) = \left\{ D_\gra \right\}$; then for every spherical root $\grs$ it holds
        \[
            c(D_\gra, \grs ) = \langle \gra^\vee , \grs \rangle.
        \]

    \item[(11)] $\scrS = (\grS, S^p, \mbA)$ is the \textit{spherical system} of $M$, where $\mbA$ has to be thought of as an abstract set together with the pairing $c : \mbA \times \grS \to \mbZ$. This is the combinatorial datum which expresses a wonderful variety: each wonderful variety is uniquely determined by its spherical system (see \cite{Lo}). There is also an abstract combinatorial definition of spherical system (see \cite{L1} or \cite{BrL}), introduced in order to obtain the classification of wonderful varieties (\textit{Luna's conjecture}): the geometrical realizability of spherical systems has been checked in many cases (see \cite{Br1}, \cite{BCF}, \cite{BrP1}, \cite{BrP2}, \cite{L1}, \cite{Pe1}) and recently a general proof which avoids a case-by-case approach has been proposed in \cite{CF}.
A very useful tool to represent spherical systems are \textit{spherical diagrams}, obtained adding information to the Dynkin diagram of $\Phi$ (see \cite{BrL}).

	\item[(12)] Denote $\gro : \Pic(M) \to \calX(B)$ and $\psi: \Pic(M) \to \calX(H)$ the restrictions to the closed and to the open orbit; then we get a commutative diagram
\[
    \xymatrix
    {
        \Pic(M)  \ar@{->}^{\psi}[r] \ar@{->}_{\gro}[d] & \calX(H) \ar@{->}[d] \\
        \calX(B)  \ar@{->}[r] & \calX(B\cap H)
    }
\]
which identifies $\Pic(M)$ with the fiber product \cite[Prop. 2.2.1]{B3}
\[
 \qquad \quad   \calX(B) \times_{\calX(B\cap H)} \calX(H) = \left\{ (\grl, \chi) \in \calX(B)\times \calX(H) \sthat \grl \ristretto_{B\cap H} = \chi\ristretto_{B\cap H} \right\}.
\]
On the combinatorial level, the map $\gro$ is described on colors as follows (\cite[Thm. 2.2]{Fo}):
\[
    \gro(D) = \left\{  \begin{array}{cl}
                                    \sum_{D\in \grD(\gra)} \gro_\gra & \textrm{if } D \in \mbA \cup \grD^b \\
                                    2 \gro_\gra & \textrm{if } D = D_\gra \in \grD^{2a}
                       \end{array}   \right.
\]
where $\gro_\gra$ is the fundamental dominant weight associated to $\gra\in S$.

    \item[(13)] If $\grS' \subset \grS$ is a subset of spherical roots, then the \textit{localization at $\grS'$} of $M$ is the $G$-stable subvariety
        \[
            M_{\grS'} = \bigcap_{\grs \in \grS \senza \grS'} M^\grs:
        \]
        it is a wonderful variety whose spherical system is $\scrS' = (\grS', S^p, \mbA')$, where $\mbA' = \bigcup_{\gra \in S \cap \grS'} \grD(\gra)$. Denote $\grD'$ the set of colors of $M_{\grS'}$; if $\gra \in S\cap \grS'$ and $\grb \in S \senza (\grS'\cup S^p)$ set $\grD'(\gra) = \{ \,'\!D^+_\gra, \,'\!D^-_\gra \}$ and $\grD'(\grb) = \{ \,D'_\grb \}$. Let $q: \Pic(M) \to \Pic(M_{\grS'})$ be the pullback map; then $\gro$ factors through $\Pic(M_{\grS'})$ and 
by its combinatorial description it follows that:
\begin{itemize}
 \item[-] if $\gra \in S \cap \grS'$ then $q(D^+_\gra)$ (resp. $q(D^-_\gra)$) is supported on $\,'\!D^+_\gra$ (resp. on $\,'\!D^-_\gra$) with multiplicity one, while it is not supported on $\,'\!D^-_\gra$ (resp. on $\,'\!D^+_\gra$);
 \item[-] if $\gra \in S \cap (\grS \senza \grS')$ then $q(D^+_\gra)$ and $q(D^-_\gra)$ are supported on $D'_\gra$ with multiplicity one;
 \item[-] if $\gra \in S\cap \frac{1}{2}\grS'$, then $q(D_\gra) = D'_\gra$;
 \item[-] if $\gra \in S \cap \frac{1}{2}(\grS \senza \grS')$, then $q(D_\gra) = 2 D'_\gra$;
 \item[-] if $\gra \in S^b$, then $q(D_\gra)$ is supported on $D'_\gra$ with multiplicity one and on at most one more color.
\end{itemize}

    \item[(14)] $M$ is said to be \textit{strict} if the stabilizer of any point is self-normalizing; equivalently, we will say also that $H$ is strict. A wonderful variety is strict if and only if it can be embedded in a simple projective space (see \cite{Pe}).

    \item[(15)] Consider the following sets of spherical roots
        \begin{align*}
            \qquad  \grS^D_\ell &= \left\{ \grs \in \grS \senza S \sthat \begin{array}{c} \textrm{there exists a rank one wonderful variety} \\ \textrm{whose spherical system is } (2\grs, S^p, \vuoto) \end{array} \right\},\\
            \qquad  \grS^S_\ell &= \left\{ \grs \in S \cap \grS \sthat c(D_\gra^+, \grs ) = c(D_\gra^-, \grs ) \; \forall \grs\in \grS \right\};
        \end{align*}
        set $\grS_\ell = \grS^D_\ell \cup \grS^S_\ell$ the \textit{set of loose spherical roots}. Loose spherical roots of the first kind are easily described, they are those of the following types (where $S = \left\{ \gra_1, \ldots, \gra_n \right\}$ and simple roots are labelled as in Bourbaki):
            \begin{itemize}
                \item[-] spherical roots $\grs = \gra_{i + 1} + \ldots + \gra_{i+r}$ with support of type $\sf{B}_r$ and with $\gra_{i+r} \in S^p$;
                \item[-] spherical roots $\grs = 2\gra_{i+1} + \gra_{i+2}$ with support of type $\sf{G}_2$.
            \end{itemize}
        For every $\grs \in \grS_\ell$, it is defined a $G$-equivariant automorphism $\grg(\grs)\in \Aut_G(M)$ of order 2 which fixes pointwise the $G$-stable divisor $M^\grs$ associated to $\grs$. If $\grs \in \grS^D_\ell$, then $\grg(\grs)$ acts trivially on $\grD$, while if $\grs \in \grS^S_\ell$, then $\grg(\grs)$ exchanges $D_\grs^+$ and $D_\grs^-$ and acts trivially on $\grD\senza \grD(\grs)$. Such automorphisms commute and generate $\Aut_G(M)$ (see \cite{Lo}).

By the natural identification $\Aut_G(M) = N_G(H)/H$, it follows that
            \begin{itemize}
                \item $H$ is self-normalizing if and only if $\grS_\ell = \vuoto$;
                \item $H$ is spherically closed if and only if $\grS^D_\ell = \vuoto$;
                \item $H$ is strict if and only if $S\cap \grS = \vuoto$ and $\grS_\ell = \vuoto$.
            \end{itemize}
        In particular, if $S \cap \grS = \vuoto$, then $H$ is self-normalizing if and only if it is spherically closed if and only if it is strict.

\item[(16)] $\grS(G)$ denotes the set of the \textit{spherical roots of} $G$: its elements are the spherical roots of all possible rank one $G$-wonderful varieties. Following the classification of the latter (see \cite{Ah}), such set is classified for any $G$ (see \cite{L1} and \cite{BrL}).

\end{itemize}

\section{Projective varieties with an open $B$-orbit.}

Let $V$ be a $G$-module and let $X \subset \mbP(V)$ be a projective variety with an open $B$-orbit; denote $G/H \incluso X$ the open orbit. Since $X$ contains finitely many $B$-orbits \cite[Cor. 2.6]{Kn3}, every $G$-orbit in $X$ is spherical. Let $p: \widetilde{X} \to X$ be the normalization; then $\wt{X}$ is a complete spherical variety with the same open orbit of $X$ whose orbits are naturally in bijection with those of $X$:

\begin{prop} [\cite{T} Prop. 1] \label{timashev}
The normalization morphism $p: \widetilde{X} \to X$ induces a bijection between $G$-orbits.
\end{prop}

If $Z\subset X$ is an orbit, set $Z' = p^{-1}(Z)\subset \wt{X}$ the corresponding orbit. Denote $Z_B\subset Z$ and $Z'_B \subset Z'$ the open $B$-orbits; fix base points $z_0 \in Z_B$ and $z'_0 \in Z'_B$ so that we have isomorphisms
\[
    Z'\simeq G/K', \qquad \qquad Z\simeq G/K
\]
with $K' \subset K$. Let's recall a result which will be useful in the following:

\begin{teo} [\cite{BP} Prop 5.1 and Cor. 5.2] \label{Brion}
Let $H$ be a spherical subgroup of $G$.
\begin{itemize}
    \item[i)] The algebraic group $N_G(H)/H$ is diagonalizable; moreover, if $H^0$ is the identity component of $H$, then $N_G(H) = N_G(H^0)$.
    \item[ii)] If $B$ is any Borel subgroup such that $BH$ is open in $G$, then $N_G(H)$ equals the right stabilizer of $BH$.
\end{itemize}
\end{teo}

Combing back to our situation, then we obtain:

\begin{cor} \label{stabilizzatori sobri}
$K'\subset K$ is a normal subgroup with finite index; in particular $K/K'$ is a finite diagonalizable group.
\end{cor}

\begin{proof}
Since $p$ is a finite morphism, it preserves dimensions of orbits: so we have $\dim(K') = \dim(K)$. Then $K'\subset K$ implies $(K')^0=K^0$ and we obtain
\[
    N_G(K') = N_G((K')^0) = N_G(K^0) = N_G(K).
\]
This shows
\[
    (K')^0 = K^0 \subset K' \subset K \subset N_G(K) = N_G(K')
\]
and the claim follows.
\end{proof}

\begin{lem} \label{indice}
Let $K' \subset K$ be two spherical subgroups of $G$ with $K'$ normal in $K$; fix a Borel subgroup $B$ such that $BK'$ is open in $G$ and consider the projection $\pi : G/K' \to G/K$. Then $\pi^{-1}(BK/K) = BK'/K'$ and $\pi^*: \grL_{G/K} \to \grL_{G/K'}$ identifies $\grL_{G/K}$ with a sublattice of $\grL_{G/K'}$ such that
\[
\Quot {\grL_{G/K'}} {\grL_{G/K}} \; \simeq \; \calX \left( \quot{K}{K'} \right).
\]
\end{lem}

\begin{proof}
First claim follows by the equality $BK' = BK$, which stems immediately from Theorem \ref{Brion} ii).

If $B'\subset B$, denote $\calX(B)^{B'}$ the kernel of the restriction $\calX(B) \to \calX(B')$, which is surjective by the following argument: if $U\subset B$ is the unipotent radical, then $\calX(B) = \calX(B/U)$ and $\calX(B') = \calX(B'/B'\cap U)$ and the restriction $\calX(B/U) \to \calX(B'/B'\cap U)$ is surjective since $B'/B'\cap U \subset B/U$ is a diagonalizable subgroup of the torus $B/U$.

By definition, we have isomorphisms $\grL_{G/K} \simeq \calX(B)^{B\cap K}$ and $\grL_{G/K'} \simeq \calX(B)^{B\cap K'}$; thus the restriction gives a surjective homomorphism
\[
\grL_{G/K'} \to \calX(B\cap K)^{B\cap K'} = \calX\left(\quot{B\cap K}{B\cap K'}\right)
\]
whose kernel is $\grL_{G/K}$. On the other hand $BK/K \simeq B/(B\cap K)$ and $BK'/K' \simeq  B/(B\cap K')$,
hence the equality $BK' = BK$ implies
\[
\Quot {B\cap K} {B\cap K'} \simeq \quot{K}{K'}.
\]
Therefore we get
\[
\Quot {\grL_{G/K'}} {\grL_{G/K}} \; \simeq \; \calX \left( \Quot {B\cap K} {B\cap K'} \right) \; \simeq \; \calX \left( \quot {K} {K'} \right).
\]
\end{proof}

\vspace{5pt}

Going back to our situation, since $K/K'$ is a finite diagonalizable group, it is isomorphic to its character group and we get the following corollary.
\begin{cor} \label{indice2}
Let $G/K \simeq Z \subset X$ be an orbit and let $G/K' \simeq Z' = p^{-1}(Z)$, with $K' \subset K$; then
\[
\Quot {\grL_{Z'}}{\grL_Z} \; \simeq \; \quot{K}{K'}.
\]
\end{cor}

Fix a closed orbit $Y\subset X$. Since parabolic subgroups are self-normalizing, $Y$ and $p^{-1}(Y)$ are isomorphic; from now on, we will denote both of them with the same letter $Y$. Let $y = [v^-] \in Y^{B^-}$ be the unique fixed point by $B^-$ (where $v^- \in V$ is a lowest weight vector) and let $\eta\in (V^*)^{(B)}$ be a highest weight vector such that $\langle \eta, v^-\rangle = 1$. If $P$ is the stabilizer of $[\eta]$, then $P$ and $\Stab(y)$ are opposite parabolic subgroups; denote $L = P \cap \Stab(y)$ the associated Levi subgroup.

Consider the affine $P$-stable open subset $X_0 = X \cap \mbP(V)_\eta$ defined by the non-vanishing of $\eta$ and recall that there exists an affine closed $L$-stable subvariety $S_X \subset X_0$ containing $y$ and possessing an open $(B\cap L)$-orbit such that the multiplication morphism
\[  \begin{array} {ccc}
                P^u \times S_X & \lra & X_0 \\
                (g,s) & \longmapsto & gs,
        \end{array}
\]
(where $P^u$ denotes the unipotent radical of $P$) is a $P$-equivariant isomorphism \cite[Prop. 1.2]{BLV}. Since $k[S_X /\!/ L] = k[S_X]^L = k$, we get that $S_X$ possesses a unique closed $L$-orbit, namely the fixed point $y$.

If $D\in \DGH$, denote $\ol{D}$ and $\wt{D}$ its closure respectively in $X$ and in $\wt{X}$.

\begin{lem} \label{chiusura colori}
Let $D\in \DGH$, let $Z\subset X$ be an orbit and set $Z' = p^{-1}(Z)$. Then $\ol{D} \supset Z$ if and only if $\wt{D}\supset Z'$.
\end{lem}

\begin{proof}
By Lemma \ref{indice}, $p^{-1}(Z_B) = Z'_B$ is the open $B$-orbit of $Z'$. Suppose that $\ol{D}\supset Z$ and fix $z_0 \in Z_B$; if $z'_0 \in p^{-1}(z_0) \cap \wt{D}$, then we obtain $Z'_B = Bz'_0 \subset \wt{D}$, which implies $Z' \subset \wt{D}$. Suppose conversely that $\wt{D}\supset Z'$: then $\ol{D}=p(\wt{D})\supset p(Z') = Z$.
\end{proof}

\begin{lem} \label{colori Xl}
Let $D\in \DGH$; then $\ol{D} \supset Y$ if and only if $\eta\ristretto_D \neq 0$.
\end{lem}

\begin{proof}
Notice that $P\ol{D}$ is closed and that it contains $Y$ if and only if $\ol{D}$ contains $Y$.
If $\eta\ristretto_D \neq 0$, then $\ol{D}\cap X_0 \neq \vuoto$, thus $P\ol{D}\cap S_X$ is non-empty, $L$-stable and closed in $S_X$: hence $y\in P\ol{D}\cap S_X$, which implies that $Y_B = By \subset P\ol{D}$. If conversely $\eta \ristretto_D = 0$, then $\ol{D}\subset \mbP(\ker(\eta))$ and we get $\ol{D} \not\supset Y$.
\end{proof}

\begin{oss}
Combining previous lemmas it follows that the set of colors $\grD_Y(\wt{X})\subset \DGH$ whose closure in $\wt{X}$ contains the closed orbit $Y$ is
\[
        \grD_Y(X) = \left\{ D\in \DGH \sthat \eta\ristretto_D \neq 0 \right\}.
\]

If $X$ possesses a unique closed orbit $Y$, then previous lemmas allow us to compute the \textit{colored cone} of the normal embedding $G/H \incluso \wt{X}$: this is the couple $\big(\calC_Y(\wt{X}), \grD_Y(\wt{X})\big)$, where $\calC_Y(\wt{X}) \subset (\grL^\vee_{G/H})_\mbQ$ is the cone generated by the elements $\rho_{\wt{X}}(D)$, where $D\subset \wt{X}$ is any $B$-stable (possibly $G$-stable) prime divisor which contains $Y$ (see \cite{Kn}). In fact, since $\wt{X}$ is simple and complete, then $\calC_Y(\wt{X})$ contains the $G$-invariant valuation cone $\calV_{G/H}$ (\cite[Thm. 5.2]{Kn}): therefore $\calC_Y(\wt{X})$ is the cone generated by $\calV_{G/H}$ together with $\rho_{X}(\grD_Y(X))$.
\end{oss}

\vspace{5pt}

Let $Z\subset X$ be an orbit. Set $Z_0 = \ol{Z} \cap \mbP(V)_\eta$, $\wt{X}_0 = p^{-1}(X_0)$ and $Z'_0 = p^{-1}(Z_0)$; then $Z_0 = X_0 \cap \ol{Z}$ and $Z'_0 = \wt{X}_0 \cap \ol{Z'}$. Considering the rings of functions, we get a commutative diagram
\[
    \xymatrix
    {
        k[X_0]  \ar@{^{(}->}[r] \ar@{->>}[d] & k[\wt{X}_0] \ar@{->>}[d] \\
        k[Z_0]  \ar@{^{(}->}[r] & k[Z'_0]
    }
\]

\begin{teo} [\cite{Kn} Thm. 2.3] \label{suriettivita semiinvarianti}
Every $B$-semiinvariant function $f\in k[Z_0]^{(B)}$ (resp. $f\in k[Z'_0]^{(B)}$) can be extended to a $B$-semiinvariant function $f' \in k[X_0]^{(B)}$ (resp. $f' \in k[\wt{X}_0]^{(B)}$).
\end{teo}

If $\grL$ is a finitely generated free $\mbZ$-module and $\grG \subset \grL$ is a submonoid, then the \textit{saturation} of $\grG$ in $\grL$ is the submonoid $\ol{\grG} = \grG_{\mbQ^+} \cap \grL$, where $\grG_{\mbQ^+}\subset \grL_\mbQ$ is the cone generated by $\grG$. If $\grO \subset \grL$ is a submonoid containing $\grG$, then we will say that $\grG$ is \textit{saturated in} $\grO$ if $\grG = \grG_{\mbQ^+} \cap \grO$.

\begin{prop} \label{reticoli saturi}
Fix an orbit $Z\subset X$ and set $Z' = p^{-1}(Z)$.
\begin{itemize}
        \item[i)] $\grL_{Z'}$ is the saturation of $\grL_Z$ in $\grL_{G/H}$.
        \item[ii)] If $Z \simeq Z'$, then $\ol{Z'} \subset \wt{X}$ is the normalization of $\ol{Z} \subset X$.
\end{itemize}
\end{prop}

\begin{proof}
Set $\grO(\wt{X}) = k[\wt{X}_0]^{(B)}/k^*$ and $\grO(\ol{Z'}) = k[Z'_0]^{(B)}/k^*$; in terms of colored cones, such monoids are described as follows (\cite[Thm. 3.5]{Kn}):
$$
	\grO(\wt{X}) = \grL_{G/H} \cap \calC^\vee_Y(\wt{X}), \qquad \qquad \grO(\ol{Z'}) = \grL_{Z'} \cap \calC^\vee_Y(\ol{Z'}).
$$
Since every $B$-semiinvariant function is uniquely determined by its weight up to a scalar factor, by previous theorem restriction gives an isomorphism of multiplicative monoids
\[
    \left\{f \in k[\wt{X}_0]^{(B)} \sthat f\ristretto_{Z'}\not \equiv 0 \right\} \stackrel{\sim}{\lra} k[Z'_0]^{(B)}:
\]
since the first one is saturated in $k[\wt{X}_0]^{(B)}$, we may then identify $\grO(\ol{Z'})$ with a saturated submonoid of $\grO(\wt{X})$. 
Hence we may as well identify $\grL_{Z'}$ with a sublattice of $\grL_{G/H}$: in fact every $B$-semiinvariant rational function on $\wt{X}$
(resp. on $Z'$) can be written as a quotient of two $B$-semiinvariant regular functions on $\wt{X}_0$ (resp. on $Z'_0$).

Since $\wt{X}$ is normal, $\grO(\wt{X}) \subset \grL_{G/H}$ is saturated; therefore, being saturated in $\grO(\wt{X})$, we see that $\grO(\ol{Z'}) \subset \grL_{G/H}$ as well is saturated.
Since the colored cone of a spherical variety is strictly convex (\cite[Thm. 4.1]{Kn}), $\calC^\vee_Y(\ol{Z'}) \subset (\grL_{Z'})_\mbQ$ has maximal dimension:
thus the equality
$$
\grL_{Z'} \cap \calC^\vee_Y(\ol{Z'}) = \grO({\ol{Z'}}) = \grL_{G/H} \cap \calC^\vee_Y(\ol{Z'})
$$
implies that $\grL_{Z'} \subset \grL_{G/H}$ is a saturated sublattice. Hence by $[\grL_{Z'} : \grL_Z] = [K:K'] < \infty$ we get i).

Finally, ii) stems from the fact that a $G$-stable subvariety of a spherical variety is normal \cite[Prop. 3.5]{BP} together with the fact that the restriction $p: \ol{Z'} \to \ol{Z}$ is finite and birational.
\end{proof}

\section{Faithful divisors.}

Let $V$ be a simple $G$-module and suppose $G/H \simeq Gx_0 \subset \mbP(V)$ is a spherical orbit. A necessary and sufficient condition so that a spherical subgroup $H$ arises in such a way has been given in \cite{BrL}:

\begin{prop}[\cite{BrL} Cor. 2.4.2] \label{sferically closed}
A spherical subgroup occurs as the stabilizer of a point in a simple projective space if and only if it is spherically closed.
\end{prop}

Let $M$ be the wonderful completion of $G/H$: then the embedding $G/H \incluso \mbP(V)$ extends to a morphism $\phi: M \to \mbP(V)$ and it determines a $B$-stable effective divisor $\grd \in \Pic(M)$ which is generated by global sections. Conversely, it is possible to start with an arbitrary spherically closed wonderful variety $M$ together with an effective $B$-stable divisor $\grd$ generated by global sections and to consider then the associated morphism $\phi_\grd : M \to \mbP(V)$, where $V = \langle G s \rangle ^*$ is the dual of the simple module generated by the canonical section $s \in \grG(M,\calO(\grd))$. In \cite{BrL}, there have been given necessary and sufficient conditions on such a divisor $\grd$ so that $\phi_\grd$ restricts to an embedding of the open orbit $G/H$; the aim of this subsection is to recall some facts which lead to such conditions.

Let $M$ be a wonderful variety with base point $x_0$ and set $H = \Stab(x_0)$. Set $\scrS = (\grS, S^p, \mbA)$ its spherical system and $\grD = \grD(G/H)$ its set of colors. Recall that a subset $\grD^* \subset \grD$ is said to be \textit{distinguished} if there exists $\grd\in \mbN_{>0}\grD^*$ such that $c(\grd, \grs) \geq 0$, for every $\grs\in \grS$.

If $H'\supset H$ is a sober subgroup and if $\phi : G/H \to G/H'$ is the projection, then the subset of colors
\[
    \grD_\phi = \{ D \in \grD \sthat \ol{\phi(D)} = G/H'\}
\]
is distinguished; conversely, if $\grD^*\subset \grD$ is a distinguished subset, then there exists a unique wonderful subgroup $H' \supset H$ with $H'/H$ connected such that $\grD^* = \grD_\phi$. This is the content of the following theorem:

\begin{teo} [\cite{Kn} Thm. 5.4, \cite{L1} Prop. 3.3.2, \cite{Bra} Thm. 3.1.1] \label{corrispondenza distinta}
$ $\\There is an inclusion-preserving bijection as follows
\[
        \begin{array} {ccc}
                    \left\{  \begin{array}{c}
                                                 \grD^* \subset \grD \textrm{ distinguished}
                                     \end{array} \right\}
                    & \longleftrightarrow   &
                    \left\{  \begin{array}{c}
                                                    H' \subset G \textrm{ wonderful } : \\ 		
                                                    H\subset H' \mand H'/H \textrm{ connected}
                                 \end{array} \right\}
        \end{array}
\]
Moreover, if $H'\supset H$ is a wonderful subgroup with $H'/H$ connected
and if $\grD^*\subset \grD$ is the corresponding distinguished subset, then
\begin{itemize}
	\item[i)] the projection $G/H \to G/H'$ identifies $\grD(G/H')$ with $\grD\senza \grD^*$;
	\item[ii)] the spherical system of the wonderful completion of $G/H'$ is
	$$\scrS/\grD^* = \big( \grS/\grD^*, S^p/\grD^*, \mbA/\grD^* \big),$$ defined as follows:
\begin{itemize}
   		 \item[-] $\grS/\grD^*$ is the set of indecomposable elements of the free semigroup
$$\quot{\mbN\grS}{\grD^*} = \left\{ \grs \in \mbN\grS \sthat c(D, \grs) = 0, \; \forall D\in \grD^* \right\};$$
  		 \item[-] $S^p/\grD^* = S^p \cup \left\{ \gra \in S \sthat \grD(\gra) \subset \grD^* \right\}$;
    	 \item[-] $\mbA/\grD^* = \bigcup_{\gra \in S\cap \grS/\grD^*} \mbA(\gra)$, and the pairing is obtained by restriction.
\end{itemize} \end{itemize}
\end{teo}

In the notations of previous theorem, the wonderful completion of $G/H'$ is denoted $M/\grD^*$ and it is called the \textit{quotient wonderful variety} of $M$ by $\grD^*$, while $\scrS/\grD^*$ is called the \textit{quotient spherical system} of $\scrS$ by $\grD^*$.

\begin{oss} \label{distinti saturi}
By \cite[Lemma 5.3 and Thm. 5.4]{Kn} together with \cite[Lemma 3.3.1]{L1},
there is an inclusion-preserving bijection between distinguished subsets $\grD^*\subset \grD$ and sober subgroups $H'\supset H$ such that $H'/H$ is connected. In \cite[Cor. 5.6.2]{L1} it was proved that, in case $G$ is of type ${\sf A}$, then such a subgroup $H'$ is necessarily wonderful; although this was claimed in general in \cite{L2}, a general proof (which stems from the classification of spherical systems) appeared only recently in \cite[Thm. 3.3.1]{Bra}.

Suppose that $H'\supset H$ is a sober subgroup such that $H'/H$ is connected and denote $\grD^*\subset \grD$ the distinguished subset of colors which map dominantly to $G/H'$. Denote $G/H'\incluso M'$ the canonical embedding and extend the projection $\phi:G/H \to G/H'$ to a morphism $\phi: M \to M'$; consider $\grL_{G/H'}$ as a sublattice of $\grL_{G/H}$. Set $N(\grD^*) = \grL_{G/H'}^\bot \subset (\grL^\vee_{G/H})_\mbQ$: it is a linear subspace which contains $\rho_{G/H}(\grD^*)$ and which intersects the valuation cone $\calV_{G/H}$ in a face; as a cone, $N(\grD^*)$ is generated by this face together with $\rho_{G/H}(\grD^*)$ \cite[Lemma 3.3.1]{L1}. While the valuation cone $\calV_{G/H'}$ is the image of $\calV_{G/H}$ under the quotient map $(\grL^\vee_{G/H})_\mbQ \to (\grL^\vee_{G/H'})_\mbQ$ by $N(\grD^*)$, the lattice $\grL_{G/H'}$ is identified with a sublattice of $\grL_{G/H}$ as follows 
\[
    \grL_{G/H'} = \grL_{G/H} \cap N(\grD^*)^\bot
\]
\cite[Lemma 5.3 and Thm. 5.4]{Kn}; as a consequence, $\grL_{G/H'}$ is saturated in $\grL_{G/H}$.

If $\grD'$ is the set of colors of $M'$ and if $M'_0 = M' \senza \bigcup_{\grD'} D$, then $\phi^{-1}(M'_0) = M \senza \bigcup_{\grD \senza \grD^*} D$; since the fibers of $\phi$ are complete and connected, it follows that $k[M'_0] = k[\phi^{-1}(M'_0)]$. Considering the $B$-semiinvariant functions, we get then the identification of semigroups
\[
    \quot{k[M'_0]^{(B)}}{k^*} \simeq - \quot{\mbN\grS}{\grD^*}.
\]
If $M'$ is smooth, then such semigroup is free; conversely, given any distinguished subset $\grD^* \subset \grD$, the semigroup $\quot{\mbN\grS}{\grD^*}$ is free \cite[Thm. 3.1.1]{Bra}, which means that $M'$ is necessarily smooth.
\end{oss}

\vspace{5pt}

Recall the restrictions to the closed and to the open orbit $\gro : \Pic(M) \to \calX(B)$ and $\psi:\Pic(M) \to \calX(H)$ and let $\grd = \sum_\grD n(\grd,D) D$ with $n(\grd,D)\geq 0$ for every $D\in \grD$. If $s\in \grG(M,\calO(\grd))$ is the canonical section, then the submodule $\langle Gs \rangle \subset \grG(M,\calO(\grd))$ generated by $s$ is identified with the simple module $V_{\gro(\grd)}$ (which contains a unique $H$-invariant line where $H$ acts by $\psi(\grd)$) and we get a morphism
\[
    \phi_\grd : M \to \mbP(V^*_{\gro(\grd)}).
\]
Define the \textit{support} of $\grd$ as
\[
    \Supp_\grD(\grd) = \left\{ D \in \grD \sthat n(\grd,D) > 0 \right\}.
\]
As a consequence of Theorem \ref{corrispondenza distinta}, we get the following corollary.

\begin{cor} \label{distinti stabilizzatori}
Let $M$ be a wonderful variety and let $\grd \in \mbN\grD$ be a divisor generated by global sections and consider the associated morphism $\phi_\grd : M \to  \mbP(V^*_{\gro(\grd)})$. Then the correspondence of Theorem \ref{corrispondenza distinta} gives an inclusion-preserving bijection as follows
\[
        \begin{array} {ccc}
                    \left\{  \begin{array}{c}
                                                 \grD^* \subset \grD \textrm{ distinguished } : \\
                                                    \grD^* \cap \Supp_\grD(\grd) = \vuoto
                                     \end{array} \right\}
                    & \longleftrightarrow   &
       	            \left\{  \begin{array}{c}
                                                    H' \subset G \textrm{ wonderful } : \\
                                                    H\subset H' \subset \Stab(\phi_\grd (x_0))\\
                                                    \mand H'/H \textrm{ connected}
                    \end{array} \right\}
        \end{array}
\]
\end{cor}

\begin{proof}
Let $H'\supset H$ be a wonderful subgroup with $H'/H$ connected and set $\grD^* \subset \grD$ the corresponding distinguished subset. If $G/H' \incluso M'$ is the wonderful embedding, then the projection $G/H \to G/H'$ extends to a morphism $M \to M'$ and pullback identifies $\Pic(M')$ with the submodule $\mbZ[\grD \senza \grD^*] \subset \mbZ\grD = \Pic(M)$. Thus the map $M\to \mbP(V^*_{\gro(\grd)})$ factors through a map $M' \to \mbP(V^*_{\gro(\grd)})$ if and only if $\Supp_\grD(\grd) \subset \grD \senza \grD^*$.
\end{proof}

\begin{dfn} [\cite{BrL}] \label{coppia fedele}
Let $M$ be a spherically closed wonderful variety and denote $\grD$ its set of colors. A divisor generated by global sections $\grd = \sum n(\grd,D)D \in \mbN\grD$ is called \textit{faithful} if it satisfies the following conditions:
\begin{description}
\item[(FD1)] Every non-empty distinguished subset of $\grD$ intersects $\Supp_\grD(\grd)$;
\item[(FD2)] If $\gra \in \grS_\ell$ is a loose spherical root, then $n(\grd, D_\gra^+) \neq n(\grd, D_\gra^-)$.
\end{description}
\end{dfn}

\begin{prop} [\cite{BrL} Prop. 2.4.3] \label{fedelta}
Let $M$ be a spherically closed wonderful variety and let $\grd\in \mbN\grD$. Then the morphism $\phi_\grd : M \to \mbP(V^*_{\gro(\grd)})$ restricts to an embedding $G/H \incluso \mbP(V^*_{\gro(\grd)})$ if and only if $\grd$ is faithful.
\end{prop}

\begin{proof}
Fix $v_0 \in \big(V^*_{\gro(\grd)}\big)^{(H)}_{-\psi(\grd)}$ a representative of the line $\phi_\grd (x_0)$ and suppose that $H = \Stab[v_0]$; then previous corollary implies \textbf{(FD1)}. Suppose by absurd that \textbf{(FD2)} fails and let $\gra \in \grS_\ell \subset S \cap \grS$ be a loose spherical root such that $n(\grd, D_\gra^+) = n(\grd, D_\gra^-)$. If $\grg(\gra) \in \Aut_G(M) = N_G(H)/H$ is the corresponding automorphism, then $\grg(\gra)$ exchanges $D_\gra^+$ and $D_\gra^-$ and fixes every other color $D\in \grD \senza \grD(\gra)$: therefore $\grg(\gra)$ fixes $\grd$. The action of $\Aut_G(M)$ on $\Pic(M) = \mbZ\grD \simeq \calX(B) \times_{\calX(B\cap H)} \calX(H)$ is defined extending by linearity the right action of $N_G(H)/H$ on $\grD$, i.e. by the action of $N_G(H)$ on $\calX(H)$. Therefore, if $g\in N_G(H)$ is a representative of $\grg(\gra)$, then $\psi(\grd)^g = \psi(\grd)$, i.e. $g$ moves the line $[v_0]$ in a line where $H$ acts by the same character: since $H$ is spherical, such a line is unique, thus $g \in H = \Stab[v_0]$ which is absurd.

Suppose conversely that $\grd$ is a faithful divisor. By \textbf{(FD1)} it follows that $\dim H = \dim \Stab[v_0]$, therefore by Theorem \ref{Brion} we get $H \subset \Stab[v_0] \subset N_G(H)$. Suppose by absurd that there exists $g\in \Stab[v_0] \senza H$. Then $\psi(\grd)^g = \psi(\grd)$ and the equivariant automorphism corresponding to the coset $gH$ fixes $\grd$: therefore by \textbf{(FD2)} we get that every color $D\in \Supp_\grD(\grd)$ is fixed by $g$. On the other hand, since $H$ is spherically closed, every element in $N_G(H)\senza H$ acts non-trivially on $\grD$. Take $\gra\in S$ such that $g$ moves $D \in \grD(\gra)$, then we get $\gra \in \grS_\ell \subset S \cap \grS$ and $\grD(\gra) = \left\{D, D\cdot g \right\}$: therefore $n(\grd, D) = n(\grd, D\cdot g) = 0$, which contradicts \textbf{(FD2)}.
\end{proof}

\begin{cor} \label{sistema sferico stabilizzatore}
In the same hypotheses of Corollary \ref{distinti stabilizzatori}, suppose moreover that every distinguished subset of $\grD$ intersects $\Supp_\grD(\grd)$ and set
\[
    \grS(\grd) = \left\{ \gra \in \grS_\ell \sthat \gra \not\in S \textrm { or } n(\grd, D_\gra^+) = n(\grd, D_\gra^-) \right\}.
\]
Then the spherical system of $\Stab(\phi_\grd(x_0))$ is $\scrS' = (\grS', S^p, \mbA')$, where
\[
    \grS' = \big(\grS \senza \grS(\grd)\big) \cup 2\grS(\grd) \qquad \textrm{ and } \qquad \mbA' = \bigcup_{\gra \in S\cap \grS'} \mbA(\gra).
\]
\end{cor}

\begin{proof}
For every loose spherical root $\grs\in \grS_\ell$, the quotient $M/\grg(\grs)$ is easily proved to be a wonderful variety, whose spherical system is $\scrS^* = (\grS^*, S^p, \mbA^*)$, where $\grS^* = (\grS \senza \left\{\grs\right\}) \cup \left\{2\grs\right\}$ and where $\mbA^* = \bigcup_{\gra \in S\cap \grS^*} \mbA(\gra)$. If $g\in N_G(H)$ is a representative of the coset corresponding to $\grg(\grs)$, then $M/\grg(\grs) = M(G/H_\grs)$, where $H_\grs$ is the subgroup generated by $H$ together with $g$. By the first part of the proof of previous theorem it follows that $H_\grs$ fixes $\phi_\grd(x_0)$, thus we get a commutative diagram
\[
    \xymatrix
    {
        M \ar@{->}[d] \ar@{->}[drr] & &  \\
        \quot{M}{\grg(\grs)} \ar@{->}[rr] & & \mbP(V^*_{\gro(\grd)})
    }
\]
Consider now the quotient variety $M/\grG_\grd$, where $\grG_\grd\subset \Aut_G(X)$ is the subgroup generated by the elements $\grg(\grs)$, with $\grs\in \grS(\grd)$: then, by previous discussion and by Proposition \ref{fedelta}, it follows that $M/\grG_\grd$ is a spherically closed wonderful variety endowed with a faithful divisor whose associated characters are the same of $\grd$.
\end{proof}

\begin{oss}
In the hypotheses of previous corollary, the assumption that every distinguished subset of colors intersects $\Supp_\grD(\grd)$ (which is equivalent to assume that $H$ and $\Stab((\phi_\grd(x_0))$ have the same dimension) involves no loss of generality: we can always reduce to that case considering, instead of $M$, the quotient wonderful variety $M/\grD(\grd)$, where $\grD(\grd)\subset \grD$ is the maximal distinguished subset which does not intersect $\Supp_\grD(\grd)$.
\end{oss}

\section{Orbits in $\Xl$ and in $\Xlt$.}

Let $M$ be a spherically closed wonderful variety with base point $x_0$ and set $H = \Stab(x_0)$; set $\scrS = (\grS, S^p, \mbA)$ its spherical system and $\grD = \grD(G/H)$ its set of colors. If $\grd \in \Pic(M)$ is a faithful divisor, set $V = V^*_{\gro(\grd)}$ and consider the morphism $\phi_\grd: M \to \mbP(V)$. Set $\Xl = \phi_\grd(M)$ and set $p:\Xlt \to \Xl$ the normalization; then we get a commutative diagram
\[
    \xymatrix
    {
        M \ar@{->>}[drr]_{\phi_\grd} \ar@{->>}[rr]^{\tilde{\phi}_\grd} & & \Xlt \ar@{->>}[d]^{p}   \\
         & & \Xl & \!\!\!\!\!\!\!\!\!\!\!\!\!\!\!\! \subset \mbP(V)
    }
\]
If $Y \subset \Xl$ is the closed orbit, then by Lemma \ref{chiusura colori} and Lemma \ref{colori Xl} we get that $\grD_Y(X_\grd) = \grD_Y(\wt{X}_\grd)$ is canonically identified with $\grD \senza \Supp_\grD(\grd)$.

If $W\subset M$ is an orbit, in the following $\grd_W\in \Pic(\ol{W})$ will denote the pullback of $\grd\in \Pic(M)$. Notice that if $\gra \in S$ then
$$\begin{array}{ll} \Supp_\grD(\grd) \cap \grD(\gra) \neq \vuoto & \iff  \gra \not\in S^p(Y) \\
					& \iff \Supp_{\grD(W)}(\grd_W) \cap \grD(W)(\gra) \neq \vuoto \quad \forall \; W \subset M  \end{array}$$
where $S^p(Y)$ denotes the set of simple roots associated to the closed orbit $Y$.

If $M$ is strict, then the variety $\Xl$ depends only on the support of $\grd$ \cite[Lemma 1]{BGMR}: if $\grd, \grd'$ are faithful divisors on $M$, then
\[
	X_\grd \simeq X_{\grd'} \iff \Supp_\grD(\grd) = \Supp_\grD(\grd').
\]
As will be shown by Lemma \ref{inclusioni S(d)}, this is not true if $M$ is not strict.

\begin{prop} \label{sistema sferico orbita Xlt}
Let $G/K \simeq Z\subset \Xl$ be an orbit and let $G/K' \simeq Z' = p^{-1}(Z)$; let $G/K_W \simeq W \subset M$ be any orbit which maps on $Z$ and choose the stabilizers so that $K_W \subset K' \subset K$. Then $K'$ is the maximal subgroup such that
$$
	K_W \subset K' \subset K \qquad \text{ and } \qquad K'/K_W \text{ is connected.}
$$
In particular, $Z \simeq Z'$ if and only if $K/K_W$ is connected.
\end{prop}

\begin{proof}
Set $K^* = K_W K^0$ the maximal subgroup of $K$ containing $K_W$ such that $K^*/K_W$ is connected. Since $K_W \subset K'$ and since $K^0 = (K')^0$, by Theorem \ref{Brion} we get that $K^* \subset K'$ is a normal subgroup; thus by Lemma \ref{indice} it follows that $K^* = K'$ if and only if $\grL_{G/K^*} = \grL_{Z'}$.

Consider the inclusions $\grL_Z \subset \grL_{Z'} \subset \grL_{W} \subset \grL_{G/H}$: since $\grL_{W}$ is saturated in $\grL_{G/H}$, Proposition \ref{reticoli saturi} shows that $\grL_{Z'}$ is the saturation of $\grL_Z$ in $\grL_{W}$. On the other hand, by Remark \ref{distinti saturi} it follows that $\grL_{G/K^*}$ is saturated in $\grL_{W}$: since $[\grL_{G/K^*} : \grL_Z] = [K:K^*] < \infty$, we get the equality $\grL_{G/K^*} = \grL_{Z'}$.
\end{proof}

\vspace{5pt}

Combining previous proposition together with \cite[Thm. 3.3.1]{Bra} and Corollary \ref{sistema sferico stabilizzatore} we get the following corollary.

\begin{cor}
Let $G/K \simeq Z\subset \Xl$ be an orbit and let $p^{-1}(Z) \simeq G/K'$ with $K'\subset K$.
Then $K'$ is a wondwerful subgroup. If moreover $M$ is strict, then $K$ is the spherical closure of $K'$.
\end{cor}

If $Z\subset \Xl$ is an orbit and if $Z'\subset \Xlt$ is the corresponding orbit, denote $\grS_Z, \grS_{Z'} \subset \mbN\grS$ the sets of spherical roots of the respective wonderful completions. By Corollary \ref{sistema sferico stabilizzatore} there exists a bijection between $\grS_Z$ and $\grS_{Z'}$, which associates to $\grg \in \grS_Z$ the unique $\grg'\in \grS_{Z'}$ which is proportional to $\grg$: more precisely, if $\grg \neq \grg'$, then $\grg = 2\grg'$.

\begin{dfn}
If $\grs \in \grS(G)$, then we say that:
\begin{itemize}
    \item[-] $\grs$ is of \textit{type} $\BrI$ if $\grs = \gra_{i+1} + \ldots + \gra_{i+r}$ has support of type ${\sf B}_r$;
    \item[-] $\grs$ is of \textit{type} $\BrII$ if $\grs = 2\gra_{i+1} + \ldots + 2\gra_{i+r}$ has support of type ${\sf B}_r$;
    \item[-] $\grs$ is of \textit{type} $\GI$ if $\grs = 2\gra_{i+1} + \gra_{i+2}$ has support of type ${\sf G}_2$;
    \item[-] $\grs$ is of \textit{type} $\GII$ if $\grs = 4\gra_{i+1} + 2\gra_{i+2}$ has support of type ${\sf G}_2$.
\end{itemize}
\end{dfn}

Consider a spherical root $\grs \in \grS(G)$ such that $2\grs \in \grS(G)$: following the explicit description of $\grS(G)$, such a root either is a simple root, or it is of type $\BrI$ or it is of type $\GI$. If $Z \subset \Xl$ is any orbit and if $Z'\subset \Xlt$ is the corresponding orbit, define $\grS(\grd_{Z'}) \subset \grS_{Z'}$ to be the subset of spherical roots which have to be doubled to get the spherical roots of $Z$.

\begin{lem} \label{orbite non isomorfe e radici sferiche cong}
An orbit $Z\subset \Xl$ is not isomorphic to its corresponding orbit $Z' \subset \Xlt$ if and only if $Z$ possesses a spherical root $\grg$ of the shape $\grg = 2\grs_1 + \ldots + 2\grs_k$, where $\grs_1, \ldots, \grs_k \in \grS$ are pairwise distinct elements (and where $\grg' = \grs_1 + \ldots + \grs_k \in \grS_{Z'}$).
\end{lem}

\begin{proof}
By Corollary \ref{sistema sferico stabilizzatore}, $Z$ and $Z'$ are not isomorphic if and only if $\grS(\grd_Z') \neq \vuoto$; suppose $\grg' \in \grS(\grd_Z')$. By Proposition \ref{sistema sferico orbita Xlt} the wonderful completion of $Z'$ is the quotient of a wonderful subvariety $M' \subset M$; therefore we can write $\grg' = a_1 \grs_1 + \ldots + a_k \grs_k$, where $\grs_1, \ldots, \grs_k$ are spherical roots of $M'$.

Since $2\grg'\in \grS(G)$, by the discussion preceeding the lemma $\grg'$ is either a simple root, or it is of type $\BrI$ or it is of type $\GI$. If $\grg'$ is a simple root or if it is of type $\BrI$ then it follows immediately that every $a_i$ is equal to one. Suppose instead that $\grg'$ is of type $\GI$; in order to show the thesis it is enough to consider the case wherein $M'$ is a wonderful variety whose spherical roots are all supported on a subset $S' = \{\gra_1, \gra_2\}\subset S$ of type $\sf{G}_2$. An easy computation shows that, if $\grS' = S'$ and if $\grD^*$ is any distinguished subset of colors of $M'$, then the quotient $M'/\grD^*$ never possesses $2\gra_1 + \gra_2$ as a spherical root. Therefore, if $\grg' = 2 \gra_1 + \gra_2$, it must be either $\grS' = \{2\gra_1 + \gra_2\}$ or $\grS' = \{\gra_1, \gra_1 + \gra_2\}$ and the claim follows. 
\end{proof}

\vspace{5pt}

As exemplified in the following sections (Example \ref{ex2} and Example \ref{ex3}), Proposition \ref{sistema sferico orbita Xlt} together with Corollary \ref{sistema sferico stabilizzatore} allow to compute explicitly the set of orbits of $\Xl$ and that of $\Xlt$ in terms of their spherical systems. This is further simplified by the following proposition, which shows that, given an orbit $Z\subset \Xl$, there exists a minimal orbit $W_Z \subset M$ mapping on $Z$. If $\grg = \sum_{\grs \in \grS} n_\grs \grs \in \grS_Z$, define
\[
	\Supp_\grS(\grg) = \{ \grs \in \grS \sthat n_\grs\neq 0 \}
\]
its support over $\grS$; define
\[
    \grS(Z) = \bigcup_{\grg \in \grS_Z} \Supp_\grS(\grg).
\]

\begin{prop} \label{orbita minimale}
Let $Z \subset \Xl$ be an orbit and let $W_Z \subset M$ the orbit whose closure has $\grS(Z)$ as set of spherical roots. Then $W_Z$ maps on $Z$ and and every other orbit which maps on $Z$ contains $W_Z$ in its closure.
\end{prop}

\begin{proof}
Let $W\subset M$ be an orbit mapping on $Z$ and let $\grS_{W}\subset \grS$ be the associated set of spherical roots. Since $\phi_\grd(W) = Z$, we get $\grS_Z \subset \mbN\grS_W$; this shows $\grS(Z)\subset \grS_{W}$, i.e. $W_Z \subset \ol{W}$. In order to prove that $\phi_\grd(W_Z) = Z$ it is enough to notice that $\grL_{\phi_\grd(W_Z)} = \grL_{W_Z}\cap \grL_Z = \grL_Z$.
\end{proof}

\begin{oss} \label{orbite massimali}
Unlike the symmetric case (see \cite{M}), in the general spherical case there does not need to exist a maximal orbit in $M$ mapping on a fixed orbit $Z\subset \Xl$: for instance this is shown by Example \ref{ex2} and by Example \ref{ex3}.
\end{oss}

\vspace{5pt}

Since $\grS(Z)$ depends only on $\grS_Z$ (or equivalently on $\grS_{Z'}$), we get the following corollaries.

\begin{cor}
Two orbits $W_1, W_2 \subset M$ map to the same orbit in $\Xl$ if and only if
\[
    \Quot{\grS_{W_1}}{\grD(\grd_{W_1})} = \Quot{\grS_{W_2}}{\grD(\grd_{W_2})},
\]
where $\grd_{W_i}$ is the pullback of $\grd$ to $\ol{W_i}$ and where $\grD(\grd_{W_i})$ is the maximal distinguished subset of colors of $W_i$ not intersecting the support of $\grd_{W_i}$.
\end{cor}

\begin{cor} \label{orbite non isomorfe}
Two orbits in $\Xl$ (resp. in $\Xlt$) have different sets of spherical roots; in particular two orbits in $\Xl$ (resp. in $\Xlt$) are never isomorphic.
\end{cor}

\begin{oss}
If $S\cap \grS = \vuoto$, then by its combinatorial description it follows that the restriction map to the closed orbit $\gro: \Pic(M) \to \calX(B)$ is injective: this means that the generic stabilizer $H$ never fixes two different lines in the same simple module. However, if $S\cap \grS \neq \vuoto$, it could happen that a simple module $\mbP(V)$ contains two different orbits both isomorphic to the open orbit $G/H$: previous corollary shows then that there does not exist any spherical orbit in $\mbP(V)$ containing both of them in its closure. For instance, this occurs in the following example.
\end{oss}

\begin{ex} \label{ex orbite sferiche doppie}
Consider $M = \mbP^1 \times \mbP^1$, which is a wonderful variety for $G = \mathrm{SL}(2)$, and fix the base point $([1,0],[0,1])$ so that the generic stabilizer is the maximal torus $T$ of diagonal matrices. Consider the simple module $V = k[x,y]_5$ formed by the homogeneous polynomials of degree 5: then $G[x^4 y]$ and $G[x^3 y^2]$ are distinct orbits in $\mbP(V)$ both isomorphic to the open orbit $G/T$.
\end{ex}

\section{Bijectivity in the strict case.}

Keeping the notations of previous section, suppose that $M$ is strict. 
Following lemma is a stronger version of Lemma \ref{orbite non isomorfe e radici sferiche cong}. 

\begin{lem} \label{orbite non isomorfe e radici sferiche}
Let $M$ be a strict wonderful variety and let $\grd$ be a faithful divisor on it; let $Z\subset \Xl$ be an orbit. Then $Z \not \simeq Z'$ if and only if there exists a spherical root $\grg \in \grS_Z$ of type $\BrII$ and a spherical root $\grs \in \Supp_\grS(\grg)$ of type $\BdueI$.
\end{lem}

\begin{proof}
By Lemma \ref{orbite non isomorfe e radici sferiche cong}, we may assume that $Z'$ possesses a spherical root $\grg$ of type $\BrI$ or of type $\GI$. Since $S \cap \grS = \vuoto$, it is uniquely determined a spherical root $\grs \in \Supp_\grS(\grg)$ which is of type ${\sf B}^\mathrm{I}_s$ (with $2\leq s \leq r$) in the first case and of type $\GI$ in the second case. Since $M$ is strict, the latter cannot happen; thus we are in the first case.

Suppose that $s > 2$ and $2\grg \in \grS_Z$; let $\grb\in S$ be the short root in the support of $\grs$. Since $M$ is strict, $\grb$ moves a color $D_\grb \in \grD$, while $s > 2$ implies $c(D_\grb, \tau) \geq 0$ for every $\tau \in \grS$: therefore $\{D_\grb\}$ is distinguished and by the faithfulness of $\grd$ we get  $D_\grb \in \Supp_\grD(\grd)$, which implies $\grb \not \in S^p(Y)$. But this is a contradiction 
since $2\grg \in \grS_Z$ implies $\grb \in S^p(Z) \subset S^p(Y)$.
\end{proof}

\vspace{5pt}

If $\grs\in \grS$ is a spherical root of type $\BdueI$, write $\grs = \gra^\sharp_\grs + \gra^\flat_\grs$, where $\gra^\sharp_\grs, \gra^\flat_\grs \in S$ are respectively the long simple root and the short simple root in the support of $\grs$. Since $M$ is strict, both $\gra^\sharp_\grs$ and $\gra^\flat_\grs$ move exactly one color; set $\grD(\gra^\sharp_\grs) = \{D^\sharp(\grs)\}$ and $\grD(\gra^\flat_\grs) = \{D^\flat(\grs)\}$.

\begin{lem} \label{D1D2}
Let $M$ be a strict wonderful variety and let $\grd$ be a faithful divisor on it; let $\grs \in \grS$ be a spherical root of type $\BdueI$.
\begin{itemize}
    \item[i)]  If $D^\flat(\grs) \in \Supp_\grD(\grd)$, then no orbit $Z \subset \Xl$ possesses a spherical root $\grg \in \grS_Z$ of type $\BrII$ with $\grs \in \Supp_\grS(\grg)$.
    \item[ii)] If $\Supp_\grD(\grd) \cap \{D^\sharp(\grs), D^\flat(\grs)\} = \{D^\sharp(\grs)\}$, then there exists an orbit $Z \subset \Xl$ such that $2\grs \in \grS_Z$; in particular $Z\not \simeq Z'$ and the normalization $p : \Xlt \to \Xl$ is not bijective.
\end{itemize}
\end{lem}

\begin{proof}
i). If $Z \subset \Xl$ possesses a spherical root $\grg$ of type $\BrII$ supported on $\grs$, then $\gra^\flat_\grs \in S^p(Z) \subset S^p(Y)$. But this is a contradiction since $D^\flat(\grs) \in \Supp_\grD(\grd)$ implies $\gra^\flat_\grs \not \in S^p(Y)$.

ii). Consider the rank one orbit $W\subset M$ whose unique spherical root is $\grs$. If $\grD(W)(\gra^\flat_\grs) = \{'\!D^\flat(\grs)\}$ and $\grD(W)(\gra^\sharp_\grs) = \{'\!D^\sharp(\grs)\}$, then 
\[
    \Supp_{\grD(W)}(\grd_W) \cap \{'\!D^\sharp(\grs), \, '\!D^\flat(\grs)\} = \{'\!D^\sharp(\grs)\}.
\]
Set $Z = \phi_\grd(W)$ and $Z' = p^{-1}(Z)$; set $\grD(\grd_W)\subset \grD(W)$ the maximal distinguished subset not intersecting the support of $\grd_W$. Since $c('\!D^\flat(\grs), \grs) = 0$ and since $'\!D^\sharp(\grs)$ is the unique color $D \in \grD(W)$ such that $c(D,\grs) > 0$, we get
\[
'\!D^\flat(\grs) \in \grD(\grd_W) = \{D \in \grD(W) \sthat c(D,\grs) = 0 \} \senza \Supp_{\grD(W)}(\grd_W)
\]
which shows $\grS_{Z'} = \{\grs\}$. On the other hand $\grD(Z')(\gra^\flat_\grs) = \vuoto$, thus $Z'$ is not spherically closed and $\grS_Z = \{2\grs\}$.
\end{proof}

\begin{cor}\label{esempi}
\begin{itemize}
    \item[i)] If $M$ is a strict wonderful symmetric variety and if $\grd$ is a faithful divisor on it, then the normalization $p:\Xlt \to \Xl$ is bijective.
    \item[ii)] Suppose that the Dynkin diagram of $G$ is simply laced. If $M$ is any strict wonderful variety for $G$ and if $\grd\in \Pic(M)$ is any faithful divisor, then the normalization $p:\Xlt \to \Xl$ is bijective.
    \item[iii)] If $D^\flat(\grs) \in \Supp_\grD(\grd)$ for every $\grs \in \grS$ of type $\BdueI$, then the normalization morphism $p : \Xlt \to \Xl$ is bijective.
\end{itemize}
\end{cor}

\begin{proof}
By the classification of symmetric varieties, we deduce that a strict wonderful symmetric variety never possesses a spherical root of type $\BdueI$. Then all of the claims above follow straightforward by previous lemmas.
\end{proof}

\vspace{5pt}

Another proof of Corollary \ref{esempi} i) was given in \cite{M}. Following examples show some cases wherein the conditions of Lemma \ref{orbite non isomorfe e radici sferiche} are fulfilled:

\begin{ex}
Consider the wonderful model variety $M$ of $\mathrm{Spin}(7)$, whose spherical system is expressed by the spherical diagram
\[
\begin{picture}(3600,1800)(-300,-900)
           \put(0,0){\usebox{\atwo}}\put(1800,0){\usebox{\btwo}}\put(3600,0){\usebox{\wcircle}}
\end{picture}
\]
Then the divisor $\grd = D_{\gra_2}$ is faithful. Consider the codimension one orbit $W\subset M$ having spherical root $\gra_2 + \gra_3$; following Proposition \ref{sistema sferico orbita Xlt} and Corollary \ref{sistema sferico stabilizzatore}, we get the following sequence of spherical diagrams
\[
\begin{picture}(23200,1800)(-300,-900)
           \put(0,0){\usebox{\dynkinatwo}}\put(0,0){\usebox{\wcircle}}\put(1800,0){\usebox{\btwo}}\put(3600,0){\usebox{\wcircle}} \put(5400,0){\line(1,0){2600}}\put(7300,0){\usebox{\toe}}\put(6300,500){\tiny $\tilde{\phi}_\grd$}
           \put(9800,0){\usebox{\dynkinatwo}}\put(9800,0){\usebox{\wcircle}}\put(11600,0){\usebox{\btwo}}
           \put(15200,0){\line(1,0){2600}}\put(17100,0){\usebox{\toe}}\put(16100,500){\tiny $p$}
           \put(19600,0){\usebox{\dynkinatwo}}\put(19600,0){\usebox{\wcircle}}\put(21400,0){\usebox{\dynkinbtwo}}\put(21400,0){\usebox{\gcircletwo}}
\end{picture}
\]
where the first one represents the orbit $W \subset M$, the second one represents the orbit $\tilde{\phi}_\grd(W) \subset \Xlt$ and the third one represents the orbit $\phi_\grd(W)\subset \Xl$.
\end{ex}

\begin{table}[t]
\caption{Example \ref{ex2} , $\grd = D_{\gra_2}$.}
\begin{tabular}{|c|c|c|c|c|}
    \hline
Maximal & Minimal & \multirow{2}{*}{Orbit in $\Xlt$} & \multirow{2}{*}{Orbit in $\Xl$} & \multirow{2}{*}{$\grS(\grd_{Z'})$} \\
Orbits & Orbit &  &  & \\
    \hline


& &
\multirow{3}{*}{\begin{picture}(7600,1800)(-200,-700)
\put(0,0){\usebox{\dynkinbfive}}\multiput(0,0)(1800,0){3}{\usebox{\atwo}}\put(5400,0){\usebox{\gcircle}}\put(7200,0){\usebox{\aprime}}
\end{picture}} &
\multirow{3}{*}{\begin{picture}(7600,1800)(-200,-700)
\put(0,0){\usebox{\dynkinbfive}}\multiput(0,0)(1800,0){3}{\usebox{\atwo}}\put(5400,0){\usebox{\gcircle}}\put(7200,0){\usebox{\aprime}}
\end{picture}} & \\
$\left\{ 1,2,3,4,5 \right\}$ & $\left\{ 1,2,3,4,5 \right\}$ & & & $\vuoto$ \\
& & & & \\
    \hline


& &
\multirow{3}{*}{\begin{picture}(7600,1800)(-200,-700)
\put(0,0){\usebox{\dynkinbfive}}\multiput(0,0)(1800,0){3}{\usebox{\atwo}}\put(5400,0){\usebox{\gcircle}}\put(7200,0){\usebox{\wcircle}}
\end{picture}} &
\multirow{3}{*}{\begin{picture}(7600,1800)(-200,-700)
\put(0,0){\usebox{\dynkinbfive}}\multiput(0,0)(1800,0){3}{\usebox{\atwo}}\put(5400,0){\usebox{\gcircle}}\put(7200,0){\usebox{\wcircle}}
\end{picture}} & \\
$\left\{ 1,2,3,4 \right\}$ & $\left\{ 1,2,3,4 \right\}$ & & & $\vuoto$ \\
& & & & \\
    \hline


& &
\multirow{3}{*}{\begin{picture}(7600,1800)(-200,-700)
        \put(0,0){\usebox{\dynkinbfive}}\multiput(0,0)(1800,0){2}{\usebox{\atwo}}\put(5400,0){\usebox{\gcircle}}
\end{picture}} &
\multirow{3}{*}{\begin{picture}(7600,1800)(-200,-700)
        \put(0,0){\usebox{\dynkinbfive}}\multiput(0,0)(1800,0){2}{\usebox{\atwo}}\put(5400,0){\usebox{\gcircletwo}}
\end{picture}} & \\
$\left\{1,2,4,5 \right\}$ & $\left\{1,2,4 \right\}$ & & & $\{\gra_4 + \gra_5\}$ \\
& & & & \\
    \hline


& &
\multirow{3}{*}{\begin{picture}(7600,1800)(-200,-700)
\put(0,0){\usebox{\dynkinbfive}}\put(1800,0){\usebox{\gcircle}}\put(5400,0){\usebox{\wcircle}}
\end{picture}} &
\multirow{3}{*}{\begin{picture}(7600,1800)(-200,-700)
\put(0,0){\usebox{\dynkinbfive}}\put(1800,0){\usebox{\gcircle}}\put(5400,0){\usebox{\wcircle}}
\end{picture}} & \\
$\left\{1,2,3,5 \right\}$ & $\left\{1,2 \right\}$ & & & $\vuoto$ \\
& & & & \\
    \hline


& &
\multirow{3}{*}{\begin{picture}(7600,1800)(-200,-450)
        \put(0,0){\usebox{\dynkinbfive}}\put(0,0){\usebox{\wcircle}}\put(1800,0){\usebox{\gcircle}}
\end{picture}} &
\multirow{3}{*}{\begin{picture}(7600,1800)(-200,-450)
        \put(0,0){\usebox{\dynkinbfive}}\put(0,0){\usebox{\wcircle}}\put(1800,0){\usebox{\gcircletwo}}
\end{picture}} & \\
$\{2,3,4,5 \}$ & $\{2,4 \}$ & & & $\left\{ \sum_{i=2}^5 \gra_i \right\}$ \\
& & & & \\
    \hline


&
\multirow{4}{*} {$\vuoto$}
&
\multirow{4}{*} {\begin{picture}(7600,1800)(-200,-700)
        \put(0,0){\usebox{\dynkinbfive}}\put(1800,0){\usebox{\wcircle}}
\end{picture}} &
\multirow{4}{*} {\begin{picture}(7600,1800)(-200,-700)
        \put(0,0){\usebox{\dynkinbfive}}\put(1800,0){\usebox{\wcircle}}
\end{picture}} &
\multirow{4}{*} {$\vuoto$} \\
$\{1,3,4,5 \}$ & & & & \\
$\{2,3,5\}$ & & & & \\
& & & & \\

    \hline

\end{tabular}
\end{table}

\begin{ex} \label{ex2}
Consider the wonderful model variety $M$ of $\mathrm{SO}(11)$, whose spherical system is expressed by the spherical diagram
\[
\begin{picture}(7200,1800)(-300,-900)
           \put(0,0){\usebox{\atwo}}
           \put(1800,0){\usebox{\atwo}}
           \put(3600,0){\usebox{\atwo}}
           \put(5400,0){\usebox{\btwo}}
           \put(7200,0){\usebox{\aprime}}
\end{picture}
\]
Then the divisor $\grd = D_{\gra_2}$ is faithful. See Table 1 for a full list of the orbits in $\Xl$ and in $\Xlt$ (for simplicity, in the table orbits in $M$ are described by giving a subset of its spherical root index set).
\end{ex}

\vspace{5pt}

As illustrated by previous examples, main examples of strict wonderful varieties possessing a faithful divisor $\grd$ such that the normalization $p:\Xlt \to \Xl$ is not bijective arise from the context of wonderful model varieties (see \cite{L2}); as will be shown in the following, the case of a general strict wonderful variety substantially follows from this special case.

Consider a strict wonderful variety $M$ and let $\grd$ be a faithful divisor on it. Let $\grs \in \grS$ be a spherical root of type $\BdueI$ and set $\grG(\grs)$ the connected component of the Dynkin diagram of $G$ where $\grs$ is supported. If $\grG(\grs)$ is of type ${\sf B}$ or ${\sf C}$, number the simple roots in $\grG(\grs)$ which are not in $S^p$ starting from the extreme of the diagram which contains the double link.

If $\{D^\flat_\grs, D^\sharp_\grs\}$ contains a distinguished subset, then by Lemma \ref{D1D2} we get that there is no orbit $Z \subset \Xl$ possessing a spherical root $\grg$ of type $\BrII$ with $\grs \in \Supp_\grS(\grg)$ if and only if $D^\flat_\grs  \in \Supp_\grD(\grd)$. For instance, this is the case if one of the following conditions is fulfilled:
\begin{itemize}
    \item[-]   $\grG(\grs)$ is of type ${\sf B}$ or ${\sf C}$ and $\grs$ is the unique spherical root supported on $\gra_2$;
    \item[-]   $\grG(\grs)$ is of type ${\sf C}$ and $2\gra_2 \in \grS$.
\end{itemize}

Suppose that $\{D^\flat_\grs, D^\sharp_\grs\}$ does not contain any distinguished subset.
If $\grG(\grs)\neq {\sf F}_4$, then there exists $\tau \in \grS$ supported on $\gra_2$ different both from $\grs$ and from $2\gra_2$; by a case-by-case check, it turns out that either $\tau$ has support of type ${\sf A}_2$ or $\grG(\grs)$ is of type ${\sf C}$ and $\tau$ has support of type ${\sf A}_1 \times {\sf A}_1$. Thus the spherical diagram of $M$ in $\grG(\grs)$ has one of the following shapes:
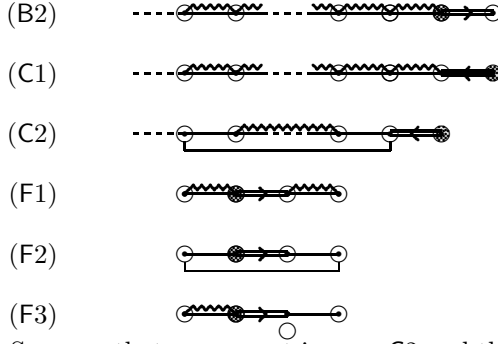
\begin{itemize}
    \item[(${\sf B}1$)] \qquad
            \begin{picture}(12600,1800)(-300,-300)
                \thicklines\multiput(0,0)(400,0){5}{\line(1,0){200}}
                \put(1800,0){\usebox{\atwoseq}}
                \put(9000,0){\usebox{\btwo}}
                \put(10800,0){\usebox{\aprime}}
            \end{picture}

    \item[(${\sf B}2$)] \qquad
            \begin{picture}(12600,1800)(-300,-300)
                \thicklines\multiput(0,0)(400,0){5}{\line(1,0){200}}
                \put(1800,0){\usebox{\atwoseq}}
                \put(9000,0){\usebox{\atwo}}
                \put(10800,0){\usebox{\btwo}}
                \put(12600,0){\usebox{\wcircle}}
            \end{picture}

    \item[(${\sf C}1$)] \qquad
        \begin{picture}(10800,1800)(-300,-300)
            \thicklines\multiput(0,0)(400,0){5}{\line(1,0){200}}
            \put(1800,0){\usebox{\atwoseq}}
            \put(9000,0){\usebox{\atwo}}
            \put(10800,0){\usebox{\leftbiedge}}
            \put(12600,0){\usebox{\gcircle}}
        \end{picture}

    \item[(${\sf C}2$)] \qquad
        \begin{picture}(10800,1800)(-300,-300)
            \thicklines\multiput(0,0)(400,0){5}{\line(1,0){200}}
            \multiput(1800,0)(7200,0){2}{\usebox{\wcircle}}
            \multiput(1800,-300)(7200,0){2}{\line(0,-1){300}}
            \put(1800,-600){\line(1,0){7200}}
            \multiput(1800,0)(5400,0){2}{\usebox{\edge}}
            \put(3600,0){\usebox{\shortam}}
            \put(9000,0){\usebox{\leftbiedge}}
            \put(10800,0){\usebox{\gcircle}}
        \end{picture}

    \item[(${\sf F}1$)] \qquad
        \begin{picture}(10800,1800)(-2100,-300)
            \put(0,0){\usebox{\dynkinf}}
            \multiput(0,0)(3600,0){2}{\usebox{\atwo}}
            \put(1800,0){\usebox{\gcircle}}
        \end{picture}

    \item[(${\sf F}2$)] \qquad
        \begin{picture}(10800,1800)(-2100,-300)
            \put(0,0){\usebox{\dynkinf}}
            \multiput(0,0)(1800,0){4}{\usebox{\wcircle}}
            \multiput(0,-300)(5400,0){2}{\line(0,-1){300}}
            \put(0,-600){\line(1,0){5400}}
            \put(1800,0){\usebox{\gcircle}}
        \end{picture}

    \item[(${\sf F}3$)] \qquad
        \begin{picture}(10800,1800)(-2100,-300)
            \put(0,0){\usebox{\dynkinf}}
            \put(0,0){\usebox{\atwo}}
            \put(1800,0){\usebox{\gcircle}}
            \put(3600,0){\usebox{\aprime}}
            \put(5400,0){\usebox{\wcircle}}
        \end{picture}

\end{itemize}

Suppose that we are not in case ${\sf C}2$ and that $\grG(\grs)$ is not of type ${\sf F}_4$: then we are substantially reduced to the case of a wonderful model variety. Let $m(\grs) \geq 3$ be the first integer such that the simple root $\gra_{m(\grs)}$ occurs in the support of one and only one spherical root with support of type ${\sf A}_2$. For $1 \leq k \leq m(\grs)$, set $\grD(\gra_k) = \{D_k\}$; set $\grD(\grs) = \{D_1, \ldots, D_{m(\grs)}\}$ and define $\grD(\grs)^{\mathrm{even}}, \grD(\grs)^{\mathrm{odd}} \subset \grD(\grs)$ as the subsets whose element index is respectively even and odd.

\begin{lem}\label{Br gen}
Let $M$ be a strict wonderful variety possessing a spherical root $\grs$ of type $\BdueI$ such that the spherical diagram of $M$ in $\grG(\grs)$ is of type ${\sf B}1$; let $\grd$ be a faithful divisor on $M$. Then there does not exist any orbit $Z\subset X_\grd$ possessing a spherical root $\grg$ of type $\BrII$ with $\grs \in \Supp_\grS(\grg)$ if and only if $D_1 \in \Supp_\grD(\grd)$ or the following conditions are both satisfied:
\begin{itemize}
    \item[i)]  $\Supp_\grD(\grd) \cap \grD(\grs)^{\mathrm{even}} = \vuoto$;
    \item[ii)] If $M$ possesses a spherical root supported on $\gra_{m(\grs)+1}$, then $m(\grs)$ is odd.
\end{itemize}
\end{lem}

\begin{proof}
By Lemma \ref{D1D2} we may assume that $\Supp_\grD(\grd) \cap \{D_1, D_2\} = \vuoto$. Notice that $\grD(\grs) \senza \{D_{m(\grs)}\}$ is distinguished and that conversely any distinguished subset which intersects $\grD(\grs)$ contains $\grD(\grs) \senza \{D_{m(\grs)}\}$.  Number the $m(\grs)$ spherical roots supported on $\{\gra_1, \ldots, \gra_{m(\grs)}\}$ from the right to the left: set $\grs_1 = 2\gra_1$ and, if $2\leq i\leq m(\grs)$, set $\grs_i = \gra_{i-1} + \gra_i$.

If $W \subset M$ is an orbit, denote $\grS' \subset \grS$ its set of spherical roots and $\grD'$ its set of colors; for $1 \leq i \leq m(\grs)$ set $\grD'(\gra_i) = \{D'_i\}$ and set $\grD'(\grs) = \{D'_1, \ldots, D'_{m(\grs)}\}$. Denote $q : \Pic(M) \to \Pic(\ol{W})$ the pullback map and observe that $q$ induces a bijection between $\grD(\grs)$ and $\grD'(\grs)$. More precisely, $q(D_i) = D'_i$ for every $1<i\leq m(\grs)$, while
\[
q(D_1) = \left\{ \begin{array}{ll} D'_1 & \textrm{if }  2\gra_1 \in \grS' \\
                                   2D'_1 & \textrm{if } 2\gra_1 \not \in \grS'\end{array} \right. :
\]
therefore, if $i \leq m(\grs)$, $\grd$ is supported on $D_i$ if and only if $\grd_W = q(\grd)$ is supported on $D'_i$.

($\Longrightarrow$) Consider the codimension one orbit $W$ whose spherical root set is $\grS' = \grS \senza \{\grs_3\}$; set $Z = \phi_\grd(W)$ and $Z' = p^{-1}(Z)$. Denote $\grD^*\subset \grD'$ the maximal distinguished subset of colors which does not intersect the support of $\grd_W$; since $D'_1 \not \in \Supp_{\grD'}(\grd_W)$ and since it is non-negative against any spherical root, we get $D'_1 \in \grD^*$.

Suppose that i) or ii) fails. Notice that, in order to show that $Z\not \simeq Z'$, it is enough to show that $D'_2 \not \in \grD^*$. On one hand, by Proposition \ref{sistema sferico orbita Xlt} together with Corollary \ref{distinti stabilizzatori} this implies $\grs \in \grL_{Z'}$: in fact $c(D',\grs)=0$ for every $D' \in \grD' \senza \{D'_2, D'_3\}$ and $D'_2 \not\in \grD^*$ implies $D'_3 \not\in \grD^*$. On the other hand, since $D'_1 \in \grD^*$, we get $\grD(Z')(\gra_1) = \grD(Z)(\gra_1) = \vuoto$, which implies that $\grs \not \in \grL_Z$. Therefore, if $D'_2 \not \in \grD^*$, then $\grs \in \grL_{Z'}\senza \grL_Z$ and $2\grs\in \grS_Z$.

Suppose first that i) fails and that $D'_2 \in \grD^*$. Then it must be either $\grD'(\grs)^{\mathrm{even}} \subset \grD^*$ or $\grD'(\grs) \senza \{D'_{m(\grs)}\} \subset \grD^*$: this follows by considering the conditions defining a distinguished subset only for $\grs_1, \grs_2, \grs_4, \ldots, \grs_{m(\grs)}$ and noticing that the minimal subsets with this property which contain $D'_2$ are
$\{D'_1\} \cup \grD'(\grs)^{\mathrm{even}}$ and, in case $m(\grs)$ is even, $\grD'(\grs) \senza \{D'_{m(\grs)}\}$. Since i) fails, the first case is not possible, while the second case is not possible because of the faithfulness of $\grd$: thus if i) fails it must be $Z\not \simeq Z'$.

Suppose now that ii) fails and that $D'_2 \in \grD^*$: thus $m(\grs)$ is even and there exists a spherical root $\grs'$ supported on $\gra_{m(\grs)+1}$. Set $m_1 := m(\grs)$ and notice that $\grs'$ has necessarily support of type ${\sf A}$. Set $m_2 > m_1+1$ the first integer such that $\gra_{m_2}$ occurs in the support of exactly one spherical root with support of type ${\sf A}$ and, proceeding similarly, define a sequence
\[
    m_1 < m_2 < \ldots < m_k 
\]
until no spherical root is supported on $\gra_{m_k+1}$. If $1 \leq j \leq m_k$, set $\grD(\gra_j) = \{D_j\}$ and $\grD'(\gra_j) = \{D'_j\}$; if $1\leq i \leq k$, set
\[
    \grD_i = \bigcup_{t = m_{i-1}+1}^{m_i} \grD(\gra_t), \quad \quad \grD'_i = \bigcup_{t = m_{i-1}+1}^{m_i} \grD'(\gra_t)
\]
(where $m_0 := 0$). Set moreover $\grD_i^{\mathrm{even}} \subset \grD_i$ and $(\grD'_i)^{\mathrm{even}} \subset \grD'_i$ the subsets whose element index $t$ is even. Define $k_0 \in \{1, \ldots, k\}$ the first integer such that $m_{k_0}$ is odd or define $k_0 = k$ otherwise. Then it is easy to show that $D'_2 \in \grD^*$ if and only if $\grD^* \cap \grD'_i = (\grD'_i)^{\mathrm{even}}$, for every $i \leq k_0$: since $\grD_{k_0}^{\mathrm{even}}\subset \grD$ is distinguished, this is impossible. Therefore if ii) fails it must be $D'_2 \not \in \grD^*$.

($\Longleftarrow$) Set $M' \subset M$ the $G$-stable prime divisor associated to the spherical root $\grs_1$ and set $W \subset M'$ the open orbit. If $Z \subset \Xl$ is an orbit possessing a spherical root $\grg$ of type $\BrII$ with $\grs \in \Supp_\grS(\grg)$, then $\grs_1 \not \in \grS(Z)$: in fact no spherical root supported on $\gra_1$ is compatible with $\grg$. Therefore by Proposition \ref{orbita minimale} such an orbit is necessarily contained in $\phi_\grd(M')$ and, in order to prove the claim, it is enough to show that it is true for any orbit which is contained in $\phi_\grd(M')$. Set $\grD^* \subset \grD'$ the maximal distinguished subset which does not intersect the support of $\grd_W$.

Suppose that both i) and ii) hold.
Then $\grD'(\grs)^{\mathrm{even}}$ is distinguished and by i) it follows that $\grD'(\grs)^{\mathrm{even}} \subset \grD^*$. Observe that $\grD^* \cap \grD'(\grs)^{\mathrm{odd}} = \vuoto$: in fact otherwise it should be $\grD'(\grs) \senza \{D'_1, D'_{m(\grs)}\} \subset \grD^*$, which contradicts the faithfulness of $\grd$. Therefore $\grD^* \cap \grD'(\grs) = \grD'(\grs)^{\mathrm{even}}$ and we get $\grs \not \in \grS(\phi_\grd(W))$: in fact, since $D'_3 \not \in \grD^*$, a spherical root $\grg\in \grS_{\phi_\grd(W)}$ with support of type ${\sf B}_r$ is necessarily a multiple of $\grs$, and this cannot happen since $c(D'_2,\grs) = 1$. To conclude, it is enough to notice that, if $Z \subset \phi_\grd(M')$ is any orbit, then $\grS(Z) \subset \grS(\phi_\grd(W))$.
\end{proof}

\begin{cor} \label{Br2 gen}
Let $M$ be a strict wonderful variety possessing a spherical root $\grs$ of type $\BdueI$ such that the spherical diagram of $M$ in $\grG(\grs)$ is of type ${\sf B}2$; let $\grd$ be a faithful divisor on $M$. Then there does not exist any orbit $Z\subset X_\grd$ possessing a spherical root $\grg$ of type $\BrII$ with $\grs \in \Supp_\grS(\grg)$ if and only if $D_1 \in \Supp_\grD(\grd)$.
\end{cor}

\begin{proof}
Let $M'$ be the wonderful variety whose spherical system is the same one of $M$ with one further spherical root $2\gra_1$: then $M$ is identified with a $G$-stable prime divisor of $M'$ and the spherical diagram of $M'$ in $\grG(\grs)$ is of the type considered in previous lemma. Denote $\grS'$ and $\grD'$ the set of spherical roots and the set of colors of $M'$; observe that the pullback map $q:\Pic(M') \to \Pic(M)$ induces an isomorphism between the sublattices generated by $\grD \senza \{D_{\gra_1}\}$ and $\grD' \senza \{D'_{\gra_1}\}$. If $D_1 \in \Supp_\grD(\grd)$ then the claim follows by Lemma \ref{D1D2}; thus we may assume $D_1 \not \in \Supp_\grD(\grd)$ and we may identify $\grd$ with a divisor $\grd'$ on $M'$ which is still faithful.

If $Z\subset \phi_{\grd'}(M')$ is an orbit possessing a spherical root $\grg$ of type $\BrII$ with $\grs \in \Supp_\grS(\grg)$, then $2\gra_1 \not \in \grS'(Z)$ and by Proposition \ref{orbita minimale} we get $Z\subset X_\grd = \phi_{\grd'}(M)$: therefore such an orbit exists in $\Xl$ if and only if it exists in $\phi_{\grd'}(M')$ and we can apply previous lemma. In order to get the claim it is enough to observe that 
if condition ii) of Lemma \ref{Br gen} holds, then (in the notations of that lemma) $\grD(\grs)^{\mathrm{even}} = q(\grD'(\grs)^{\mathrm{even}}) \subset \grD$ is distinguished: thus $\Supp_\grD(\grd) \cap \grD(\grs)^{\mathrm{even}} \neq \vuoto$ and consequently i) fails.
\end{proof}

\vspace{5pt}

If they are defined, set
\[
            e_\grs(\grd) = \min\{k \leq m(\grs) \sthat D_k \in \Supp_\grD(\grd) \cap \grD(\grs)^{\mathrm{even}} \},
\]
\[
            o_\grs(\grd) = \min\{k \leq m(\grs) \sthat D_k \in \Supp_\grD(\grd) \cap \grD(\grs)^{\mathrm{odd}} \}.
\]

\begin{lem}\label{Cr gen}
Let $M$ be a strict wonderful variety possessing a spherical root $\grs$ of type $\BdueI$ such that the spherical diagram of $M$ in $\grG(\grs)$ is of type ${\sf C}1$; let $\grd$ be a faithful divisor on $M$. Then there does not exist any orbit $Z\subset X_\grd$ possessing $2\grs$ as a spherical root if and only if the following conditions are both satisfied:
\begin{itemize}
    \item[i)] $\Supp_\grD(\grd) \cap \grD(\grs)^{\mathrm{even}} \neq \vuoto$;
    \item[ii)] If $\Supp_\grD(\grd) \cap \grD(\grs)^{\mathrm{odd}} \neq \vuoto$, then $o_\grs(\grd) \geq e_\grs(\grd)-1$.
\end{itemize}
\end{lem}

\begin{proof}
Notice that if $m(\grs)$ is even then $\grD(\grs)^{\mathrm{odd}}$ is distinguished, while if $m(\grs)$ is odd then $\grD(\grs)^{\mathrm{even}}$ is distinguished: thus at least one between $e_\grs(\grd)$ and $o_\grs(\grd)$ is well defined. By Lemma \ref{D1D2}, we may suppose $\min\{e_\grs(\grd), o_\grs(\grd)\} > 2$. Number the $m(\grs)-1$ spherical roots supported on $\{\gra_1, \ldots, \gra_{m(\grs)}\}$ from the right to left: if $i< m(\grs)$, set $\grs_i = \gra_i + \gra_{i+1}$.

If $W \subset M$ is an orbit , denote $\grS' \subset \grS$ its set of spherical roots and $\grD'$ its set of colors; for $1 \leq i \leq m(\grs)$ set $\grD'(\gra_i) = \{D'_i\}$ and set $\grD'(\grs) = \{D'_1, \ldots, D'_{m(\grs)}\}$. Denote $q : \Pic(M) \to \Pic(\ol{W})$ the pullback map and observe that $q$ induces a bijection between $\grD(\grs)$ and $\grD'(\grs)$. Since $q(D_i) = D'_i$ for every $i \leq m(\grs)$, $\grd$ is supported on $D_i$ if and only if $\grd_W = q(\grd)$ is supported on $D'_i$.

($\Longrightarrow$) Suppose that $o_\grs(\grd)$ is defined and, in case $e_\grs(\grd)$ is defined too, suppose that $o_\grs(\grd) < e_\grs(\grd) - 1$:
this implies $o_\grs(\grd)<m(\grs)$, since otherwise $\grD(\grs)^{\mathrm{even}}$ would be distinguished and it would be $e_\grs(\grd)<o_\grs(\grd)$. Consider the orbit $W \subset M$ whose spherical roots are $\grs_1, \ldots, \grs_{j}$; set $Z = \phi_\grd(W)$ and $Z' = p^{-1}(Z)$. Then the maximal distinguished subset of $\grD'$ which does not intersect the support of $\grd_W$ is
\[
    \grD^* = \grD' \senza \big(\grD'(\grs)^{\mathrm{odd}}_{\leq j+2} \cup \Supp_{\grD'}(\grd_W) \big),
\]
which by hypothesis contains $\grD'(\grs)^{\mathrm{even}}_{\leq j+1}$ (where the notations are the obvious ones); thus $\grD^* \cap \{D'_1, D'_2, D'_3\} = \{D'_2\}$. Since $c(D', \grs) = 0$ for every $D' \in \grD' \senza \{D'_1, D'_3\}$, by Proposition \ref{sistema sferico orbita Xlt} together with Corollary \ref{distinti stabilizzatori} we get $\grs \in \grL_{Z'}$. On the other hand, $D'_2 \in \grD^*$ implies $\grD(Z)(\gra_2) = \vuoto$: since $Z$ is spherically closed, we get then $\grs \not \in \grS_Z$ and $2\grs \in \grS_Z$.

($\Longleftarrow$) Suppose that $e_\grs(\grd)$ is defined and, in case $o_\grs(\grd)$ is defined too, suppose that $o_\grs(\grd)\geq e_\grs(\grd)-1$. Fix an orbit $W \subset M$, set $Z = \phi_\grd(W)$ and $Z' = p^{-1}(Z)$. We may assume that $\grs \in \grS'$, since otherwise there is nothing to prove. Set $\grD^* \subset \grD'$ the maximal distinguished subset which does not intersect the support of $\grd_W$ and notice that $2\grs \in \grS_Z$ if and only if $\grD^* \cap \{D'_1, D'_2, D'_3 \} = \{ D'_2 \}$. Such condition does not hold if $\grs_2 \not \in \grS'$ or if $\grs_3 \not \in \grS'$, since then it would be $D'_1 \in \grD^*$; thus we may assume that $\grS' \supset \{\grs_1, \grs_2, \grs_3 \}$.

Set $k < m(\grs)$ the maximum such that $\grs_i \in \grS'$ for every $i \leq k$. By considering the conditions defining a distinguished set only for $\grs_1, \ldots, \grs_k$ it follows that, if $D'_2 \in \grD^*$, then either $\grD'(\grs)_{\leq k} \subset \grD^*$ or $\grD'(\grs)^{\mathrm{even}}_{\leq k+1} \subset \grD^*$. If we are in the first case, then we are done; suppose we are in the second case. Then it must be $e_\grs(\grd)>k+1$ and, by the hypothesis, we get $o_\grs(\grd) > k$. Since it is distinguished and it does not intersect the support of $\grd_W$, we get then $\grD'(\grs)_{\leq k}\subset \grD^*$: therefore the condition $\grD^* \cap \{D'_1, D'_2, D'_3 \} = \{ D'_2 \}$ is not satisfied whenever conditions i) and ii) hold and the claim follows.
\end{proof}

\vspace{5pt}

Combining together Lemma \ref{Br gen}, Corollary \ref{Br2 gen} and Lemma \ref{Cr gen}, we get the following theorem (the cases wherein the spherical diagram of $M$ in $\grG(\grs)$ is of type ${\sf C}2$, ${\sf F}1$, ${\sf F}2$ or ${\sf F}3$ are easily treated directly).

\begin{teo} \label{teorema caso stretto}
Let $M$ be a strict wonderful variety and let $\grd$ be a faithful divisor on it. Then the normalization $p: \Xlt \to \Xl$ is bijective if and only if the following conditions are fulfilled, for every spherical root $\grs\in \grS$ of type $\BdueI$:
\begin{itemize}
    \item[i)] If the spherical diagram of $M$ in $\grG(\grs)$ is of type ${\sf B}1$, then $D^\flat(\grs) \in \Supp_\grD(\grd)$ or the following conditions are both satisfied:
                \begin{itemize}
                        \item  $\Supp_\grD(\grd) \cap \grD(\grs)^{\mathrm{even}} = \vuoto$;
                        \item If $M$ possesses a spherical root supported on $\gra_{m(\grs)+1}$, then $m(\grs)$ is odd.
                \end{itemize}
    \item[ii)] If the spherical diagram of $M$ in $\grG(\grs)$ is of type ${\sf B}2$, then $D^\flat(\grs) \in \Supp_\grD(\grd)$.
    \item[iii)] If the spherical diagram of $M$ in $\grG(\grs)$ is of type ${\sf C}1$, then the following conditions are both satisfied
                \begin{itemize}
                        \item $\Supp_\grD(\grd) \cap \grD(\grs)^{\mathrm{even}} \neq \vuoto$;
                        \item If $\Supp_\grD(\grd) \cap \grD(\grs)^{\mathrm{odd}} \neq \vuoto$, then $o_\grs(\grd) \geq e_\grs(\grd)-1$.
                \end{itemize}
    \item[iv)] Otherwise, if $D^\sharp(\grs)\in \Supp_\grD(\grd)$, then $D^\flat(\grs)\in \Supp_\grD(\grd)$ as well.
\end{itemize}
\end{teo}

\section{Bijectivity in the non-strict case.}

Keeping the notations of previous sections, suppose that $M$ is not strict and let $\grd = \sum_\grD n(\grd, D) D$ be a faithful divisor on $M$. Suppose that $Z\subset \Xl$ is an orbit such that $\grS(\grd_Z)$ contains a non-simple spherical root $\grg$. Following examples show that, unlike from the strict case (Lemma \ref{orbite non isomorfe e radici sferiche}), $\grg$ may be as well of type $\GI$ and, in case $\grg$ is of type $\BrI$, then it does not necessarily come from a spherical root of type $\BdueI$.

\begin{ex}
Consider the wonderful variety $M$ whose spherical system is expressed by the spherical diagram
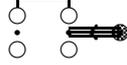
\[
\begin{picture}(3600,1800)(-300,-900)
           \put(0,0){\usebox{\aone}}
           \multiput(0,900)(1800,0){2}{\line(0,1){300}}\put(0,1200){\line(1,0){1800}}
           \put(1800,0){\usebox{\dynking}}\put(1800,0){\usebox{\aone}}\put(3600,0){\usebox{\gcircle}}
\end{picture}
\]
Then the divisor $\grd = D^+_{\gra_1}$ is faithful. Consider the codimension one orbit $W\subset M$ whose spherical roots are $\gra_2$ and $\gra_2 + \gra_3$; following Proposition \ref{sistema sferico orbita Xlt} and Corollary \ref{sistema sferico stabilizzatore}, we get the sequence of spherical diagrams
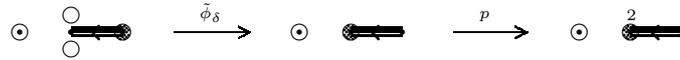
\[
\begin{picture}(23200,1800)(-300,-900)
           \put(0,0){\usebox{\wcircle}}\put(0,0){\usebox{\vertex}}
           \put(1800,0){\usebox{\dynking}}\put(1800,0){\usebox{\aone}}\put(3600,0){\usebox{\gcircle}}

            \put(5400,0){\line(1,0){2600}}\put(7300,0){\usebox{\toe}}\put(6300,500){\tiny $\tilde{\phi}_\grd$}

           \put(9800,0){\usebox{\wcircle}}\put(9800,0){\usebox{\vertex}}
           \put(11600,0){\usebox{\dynking}}\put(11600,0){\usebox{\gcircle}}

           \put(15200,0){\line(1,0){2600}}\put(17000,0){\usebox{\toe}}\put(16100,500){\tiny $p$}

           \put(19600,0){\usebox{\wcircle}}\put(19600,0){\usebox{\vertex}}
           \put(21400,0){\usebox{\dynking}}\put(21400,0){\usebox{\gcircletwo}}

\end{picture}
\]
where the first one represents the orbit $W \subset M$, the second one the orbit $\tilde{\phi}_\grd(W) \subset \Xlt$ and the third one the orbit $\phi_\grd(W)\subset \Xl$.
\end{ex}

\begin{ex} \label{ex3}
Consider the wonderful variety $M$ whose spherical system is expressed by the spherical diagram
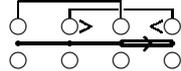
\[
	\begin{picture}(6000,2700)(-300,-900)
	\multiput(0,0)(1800,0){2}{\usebox{\edge}}\put(3600,0){\usebox{\dynkinbtwo}}\multiput(0,0)(1800,0){4}			{\usebox{\aone}} \multiput(0,900)(3600,0){2}{\line(0,1){600}}\put(0,1500){\line(1,0){3600}}
	\multiput(1800,900)(3600,0){2}{\line(0,1){300}}\multiput(1800,1200)(1900,0){2}{\line(1,0)						{1700}}\put(1800,600){\usebox{\toe}}\put(5400,600){\usebox{\tow}}
	\end{picture}
\]
Then the divisor $\grd = D^+_{\gra_1}$ is faithful. See Table 2 for a full list of the orbits in $\Xlt$ and in $\Xl$ (for simplicity, in the table orbits in $M$ are described by giving a subset of its spherical root index set).
\end{ex}

\begin{table}[t]
\caption{Example \ref{ex3}, $\grd = D^+_{\gra_1}$.}
\begin{tabular}{|c|c|c|c|c|}
    \hline
Maximal & Minimal & \multirow{2}{*}{Orbit in $\Xlt$} & \multirow{2}{*}{Orbit in $\Xl$} & \multirow{2}{*}{$\grS(\grd_{Z'})$} \\
Orbits & Orbit &  &  & \\
    \hline


& &
\multirow{3}{*}{\begin{picture}(6800,1800)(-700,-450)
\put(0,0){\usebox{\dynkinbfour}}\multiput(0,0)(1800,0){4}{\usebox{\aone}}\multiput(0,900)(3600,0){2}{\line(0,1){600}}\put(0,1500){\line(1,0){3600}}
\multiput(1800,900)(3600,0){2}{\line(0,1){300}}\multiput(1800,1200)(1900,0){2}{\line(1,0){1700}}\put(1800,600){\usebox{\toe}}\put(5400,600){\usebox{\tow}}
\end{picture}} &
\multirow{3}{*}{\begin{picture}(6800,1800)(-700,-450)
\put(0,0){\usebox{\dynkinbfour}}\multiput(0,0)(1800,0){4}{\usebox{\aone}}\multiput(0,900)(3600,0){2}{\line(0,1){600}}\put(0,1500){\line(1,0){3600}}
\multiput(1800,900)(3600,0){2}{\line(0,1){300}}\multiput(1800,1200)(1900,0){2}{\line(1,0){1700}}\put(1800,600){\usebox{\toe}}\put(5400,600){\usebox{\tow}}
\end{picture}} & \\
$\left\{ 1,2,3,4 \right\}$ & $\left\{ 1,2,3,4 \right\}$ & & & $\vuoto$ \\
& & & & \\
    \hline


& &
\multirow{3}{*}{\begin{picture}(6800,2400)(-700,-900)
\put(0,0){\usebox{\dynkinbfour}}\multiput(0,0)(1800,0){3}{\usebox{\aone}}\multiput(0,900)(3600,0){2}{\line(0,1){300}}
\put(0,1200){\line(1,0){3600}}\put(1800,600){\usebox{\toe}}\put(5400,0){\usebox{\wcircle}}
\end{picture}} &
\multirow{3}{*}{\begin{picture}(6800,2400)(-700,-900)
\put(0,0){\usebox{\dynkinbfour}}\multiput(0,0)(1800,0){3}{\usebox{\aone}}\multiput(0,900)(3600,0){2}{\line(0,1){300}}
\put(0,1200){\line(1,0){3600}}\put(1800,600){\usebox{\toe}}\put(5400,0){\usebox{\wcircle}}
\end{picture}} & \\
$\left\{ 1,2,3 \right\}$ & $\left\{ 1,2,3 \right\}$ & & & $\vuoto$ \\
& & & & \\
    \hline


& &
\multirow{3}{*}{\begin{picture}(6800,2400)(-700,-900)
\put(0,0){\usebox{\dynkinbfour}}\multiput(0,0)(3600,0){2}{\usebox{\aone}}\multiput(0,900)(3600,0){2}{\line(0,1){300}}\put(1800,0){\usebox{\wcircle}}
\put(0,1200){\line(1,0){3600}}\put(5400,600){\usebox{\tow}}\put(5400,0){\usebox{\aone}}
\end{picture}} &
\multirow{3}{*}{\begin{picture}(6800,2400)(-700,-900)
\put(0,0){\usebox{\dynkinbfour}}\multiput(0,0)(3600,0){2}{\usebox{\aone}}\multiput(0,900)(3600,0){2}{\line(0,1){300}}\put(1800,0){\usebox{\wcircle}}
\put(0,1200){\line(1,0){3600}}\put(5400,0){\usebox{\aprime}}
\end{picture}} & \\
$\left\{1,3,4 \right\}$ & $\left\{1,3,4 \right\}$ & & & $\{\gra_4\}$ \\
& & & & \\
    \hline


& &
\multirow{3}{*} {\begin{picture}(6800,1800)(-700,-900)
           \put(0,0){\usebox{\dynkinbfour}}\multiput(0,0)(5400,0){2}{\usebox{\wcircle}}\put(1800,0){\usebox{\atwo}}\put(3600,0){\usebox{\gcircle}}
\end{picture}} &
\multirow{3}{*} {\begin{picture}(6800,1800)(-700,-900)
           \put(0,0){\usebox{\dynkinbfour}}\multiput(0,0)(5400,0){2}{\usebox{\wcircle}}\put(1800,0){\usebox{\atwo}}\put(3600,0){\usebox{\gcircle}}
\end{picture}} & \\
$\left\{2,3,4 \right\}$ & $\left\{2,3,4 \right\}$ & & & $\vuoto$ \\
& & & & \\
    \hline


& &
\multirow{3}{*} {\begin{picture}(6800,1800)(-700,-900)
           \put(0,0){\usebox{\dynkinbfour}}\multiput(0,0)(1800,0){4}{\usebox{\wcircle}}
           \multiput(0,-600)(3600,0){2}{\line(0,1){300}}\put(0,-600){\line(1,0){3600}}
\end{picture}} &
\multirow{3}{*} {\begin{picture}(6800,1800)(-700,-900)
           \put(0,0){\usebox{\dynkinbfour}}\multiput(0,0)(1800,0){4}{\usebox{\wcircle}}
           \multiput(0,-600)(3600,0){2}{\line(0,1){300}}\put(0,-600){\line(1,0){3600}}
\end{picture}} & \\
$\{1, 3 \}$ & $\{1, 3 \}$ & & & $\vuoto$ \\
& & & & \\
    \hline


& &
\multirow{3}{*} {\begin{picture}(6800,1800)(-700,-900)
           \put(0,0){\usebox{\dynkinbfour}}\put(0,0){\usebox{\wcircle}}\put(3600,0){\usebox{\gcircle}}
\end{picture}} &
\multirow{3}{*} {\begin{picture}(6800,1800)(-700,-900)
           \put(0,0){\usebox{\dynkinbfour}}\put(0,0){\usebox{\wcircle}}\put(3600,0){\usebox{\gcircletwo}}
\end{picture}} & \\
$\{3,4 \}$ & $\{3,4 \}$ & & & $\{\gra_3 + \gra_4\}$ \\
& & & & \\
    \hline


&
\multirow{4}{*} {$\vuoto$}
&
\multirow{4}{*} {\begin{picture}(6800,1800)(-700,-900)
           \put(0,0){\usebox{\dynkinbfour}}\multiput(0,0)(3600,0){2}{\usebox{\wcircle}}
\end{picture}} &
\multirow{4}{*} {\begin{picture}(6800,1800)(-700,-900)
           \put(0,0){\usebox{\dynkinbfour}}\multiput(0,0)(3600,0){2}{\usebox{\wcircle}}
\end{picture}} &
\multirow{4}{*} {$\vuoto$} \\
$\{1,2,4 \}$ & & & & \\
$\{2,3\}$ & & & & \\
& & & & \\

    \hline

\end{tabular}
\end{table}

\begin{lem}\label{inclusioni S(d)} Suppose that $M$ is a spherically closed wonderful variety and let $\grd = \sum_\grD n(\grd, D) D$ be a faithful divisor on it; let $\gra \in S \cap \grS$.
\begin{itemize}
    \item [i)] If $Z\subset \Xl$ is an orbit such that $2\gra \in \grS_Z$, then $n(\grd, D^+_\gra) = n(\grd, D^-_\gra)$.
    \item [ii)] If $n(\grd, D^+_\gra) = n(\grd, D^-_\gra)$ is non-zero, then there exists an orbit $Z\subset \Xl$ such that $2\gra \in \grS_Z$.
\end{itemize}
\end{lem}

\begin{proof}
Suppose that $W\subset M$ is an orbit with set of spherical roots $\grS' \subset \grS$ and set of colors $\grD'$. If $\gra \in S \cap \grS'$, set $\grD'(\gra) = \{\,'\!D^+_\gra , \,'\!D^-_\gra \}$;
then by the description of the pullback map $q: \Pic(M) \to \Pic(\ol{W})$ it follows that
\[
n \big(\grd_W, \,'\!D^+_\gra \big) = n(\grd, D^+_\gra), \qquad n \big(\grd_W, \,'\!D^-_\gra \big) = n(\grd, D^-_\gra),
\]
where $\grd_W = q(\grd)$.

i). Let $Z\subset \Xl$ be an orbit possessing $2\gra$ as a spherical root; let $Z'=p^{-1}(Z)$ and let $W\subset M$ be an orbit which maps on $Z$. Then by Proposition \ref{reticoli saturi} we get that $\gra\in \grS_{Z'}$, while by Corollary \ref{sistema sferico stabilizzatore} together with Theorem \ref{corrispondenza distinta} we get $n(\grd_W, \,'\!D^+_\gra) = n(\grd_W, \,'\!D^-_\gra)$; by the remark at the beginning of the proof this implies the thesis.

ii). Consider the rank one orbit $W$ whose unique spherical root is $\gra$, set $Z = \phi_\grd(W)$ and $Z' = p^{-1}(Z)$. Then $\gra\in \grS_{Z'}$ is a loose spherical root and by the remark at the beginning of the proof we get $n(\grd_{Z'}, \,'\!D^+_\gra) = n(\grd_{Z'}, \,'\!D^-_\gra)$, where $\grd_{Z'}$ is the pullback of a hyperplane section and where $\grD(Z')$ is identified with a subset of $\grD(W)$.
Then by Corollary \ref{sistema sferico stabilizzatore} we get that $2\gra\in \grS_Z$.
\end{proof}

\vspace{5pt}

Suppose that $\gra \in S \cap \grS$. As shown by Example \ref{ex3}, if $n(\grd, D^+_\gra) = n(\grd, D^-_\gra) = 0$, then it may not exist any orbit $Z\subset \Xl$ possessing $2\gra$ as a spherical root; conversely, if there exists such an orbit, it may be as well $n(\grd, D^+_\gra) = n(\grd, D^-_\gra) = 0$.

As a corollary of previous lemma, we get the following sufficient conditions.

\begin{cor} \label{condizioni non normalita}
Suppose that $M$ is a spherically closed wonderful variety and let $\grd = \sum_\grD n(\grd, D) D$ be a faithful divisor on it.
\begin{itemize}
    \item[i)] If there exists $\gra \in S \cap \grS$ such that $n(\grd, D_\gra^+) = n(\grd, D_\gra^-)$ is non-zero, then the normalization $p: \Xlt \to \Xl$ is not bijective.
    \item[ii)] If the Dynkin diagram of $G$ is simply laced and if $n(\grd, D_\gra^+) \neq n(\grd, D_\gra^-)$ for every $\gra \in S \cap \grS$, then the normalization $p: \Xlt \to \Xl$ is bijective.
\end{itemize}
\end{cor}

Reasoning as in Lemma \ref{D1D2} and in Corollary \ref{esempi}, other sufficient conditions of bijectivity can be obtained imposing further conditions on the support of $\grd$ on the multiple links of the Dynkin diagram of $G$ and on the simple spherical roots of $M$.


\end{document}